\documentclass [12pt]{article}
\usepackage{amssymb,amsmath, amsthm, amscd,ifthen}
\usepackage[T2A]{fontenc}
\usepackage[cp1251]{inputenc}
\usepackage[english,russian]{babel}
\usepackage[dvips]{graphicx}
\usepackage{indentfirst}

\newtheorem{theorem}{Theorem}
\newtheorem{claim}{Claim}
\newtheorem{lemma}{Lemma}
\newtheorem{corollary}{Corollary}
\newenvironment{proof1}{\noindent{\it Proof.\,}}{\hfill$\Box$}

\theoremstyle{remark}

\newtheorem{note}{Remark}

\begin{document}
\title{On $r$-colorability of random hypergraphs}
\author{Andrey Kupavskii\footnote{This work was partially supported by the grant 09-01-00294 of
Russian Foundation for Basic Research and
the grant MD-8390.2010.1 of the Russian President.}, Dmitry A. Shabanov\footnote{This work was partially supported by Russian Foundation of Fundamental Research (grant no. 09-01-00294), by the program "`Leading Scientific Schools"' (grant no. NSh-8784.2010.1) and by the grant of President of Russian Federation (no. MK-3429.2010.1).}}
\date{}

\maketitle

\centerline{\bf\Large Abstract}

\vspace{4mm}
The work deals with the threshold for $r$-colorability in the binomial model $H(n,k,p)$ of a random hypergraph. We prove that if, for some constant $\delta\in(0,1)$,
$$
  k^{\varphi(k)}\ln n\ll r^{k-1}\le n^{(1-\delta)/2}\mbox{ and }p\le\,r^{k-1}k^{-1-\varphi(k)}\,\frac {n}{{n\choose k}},
$$
where $\varphi(k)$ is some function satisfying the relation $\varphi(k)=\Theta\left(\sqrt{\frac{\ln\ln k}{\ln k}}\right)$, then
$$
  {\sf P}\left(H(n,k,p)\mbox{ is $r$-colorable}\right)\to 1\mbox{ as }n\to\infty.
$$
This result improves the previously known results in the wide range of the parameters $r=r(n)$, $k=k(n)$.

\vspace{2mm}
\textbf{Keywords:} \emph{random hypergraph}, \emph{colorings of hypergraphs}, \emph{sparse hypergraphs}, \emph{random recoloring method}.

\vspace{4mm}

\section{Introduction and history of the problem}

The work deals with a problem concerning threshold for $r$-colorability in the binomial model of a random hypergraph. First of all, we recall the main definitions from the hypergraph theory.

\subsection{Main definitions}

\emph{A hypergraph} $H$ is a pair $H=(V,E)$, where $V=V(H)$ is some finite set (called \emph{the vertex set} of the hypergraph) and $=E(H)$ is a collection of subsets of $V$, which are called \emph{the edges} of the hypergraph. If $E\subseteq {V\choose k}$, i.e. every edge contains exactly $k$ vertices, then $H$ is called $k$-\emph{uniform}. By $K_n^{(k)}$ we denote a complete $k$-uniform hypergraph on $n$ vertices.

A vertex coloring of the vertex set $V$ is called \emph{proper} for hypergraph $H=(V,E)$ if in this coloring there is no monochromatic edges in $E$. A hypergraph $H$ is called $r$-colorable, if there exists a proper coloring with $r$ colors ($r$-\emph{coloring}) for $H$. \emph{The chromatic number} $\chi(H)$ of a hypergraph $H$ is the minimum $r$ such that $H$ is $r$-colorable.

Let $v$ be a vertex of a hypergraph $H$. \emph{The degree} of $v$ in $H$ is the number of edges of $H$ containing $v$. By $\Delta(H)$ we denote the maximum vertex degree of the hypergraph $H$. A hypergraph $H=(V,E)$ is called \textit{$l$-simple}, if every two of its distinct edges do not share more than $l$ common vertices, i.e.
$$
	\forall\; e,f\in E,\; f\ne e:\; \vert e\cap f\vert\le l.
$$
A 1-simple hypergraph is usually called \emph{simple} hypergraph. A \textit{cycle of length} 3 (3-\emph{cycle} or \emph{triangle}) in the hypergraph $H$ is a unordered set of three distinct edges $(e,f,h)$ such that $(e\cap f)\setminus h\ne\emptyset$, $(e\cap h)\setminus f\ne\emptyset$, $(h\cap f)\setminus e\ne\emptyset$.

\bigskip

In this article we study the binomial model of a random hypergraph. Given integers $n>k\ge 2$ and a real number $p\in[0,1]$,  \textit{random hypergraph $H(n,k,p)$} is a random spanning subhypergraph of the complete $k$-uniform hypergraph $K_n^{(k)}$ with the following distribution: for any spanning subhypergraph $F$ of $K_n^{(k)}$,
$$
  {\sf P}\left(H(n,k,p)=F\right)=p^{|E(F)|}(1-p)^{{n\choose k}-|E(F)|}.
$$
This definition immediately implies that every edge of $K_n^{(k)}$ is included in $H(n,k,p)$ independently with equal probability $p$.

Suppose ${\cal Q}_n$ is a property of $k$-uniform hypergraphs on $n$ vertices. We say that ${\cal Q}_n$ is an \emph{increasing property} if, for any two $k$-uniform hypergraphs $H$ and $G$ on $n$ vertices, $H$ has property ${\cal Q}_n$ and $E(H)\subseteq E(G)$ imply that $G$ has property ${\cal Q}_n$. For given function $k=k(n)\ge 2$, the function $p^*=p^*(n)$ is said to be a \emph{threshold} (or a \emph{threshold probability}) for an increasing property ${\cal Q}_n$, if
\begin{itemize}
  \item for any $p=p(n)$ such that $p\ll p^*$, ${\sf P}\left(H(n,k,p)\mbox{ has property }{\cal Q}_n\right)\to 0\mbox{ as }n\to\infty;$
  \item for any $p=p(n)$ such that $p\gg p^*$, ${\sf P}\left(H(n,k,p)\mbox{ has property }{\cal Q}_n\right)\to 1\mbox{ as }n\to\infty.$
\end{itemize}
It follows from general results of B. Bollob\'as and A. Thomason (see \cite{BolTom}) concerning monotone properties of random subsets that for any function $k=k(n)\ge 2$ and any increasing property ${\cal Q}_n$, the threshold probability exists. In this article we are concentrated on the estimating the threshold probability for $r$-colorability of $H(n,k,p)$, i.e. for an increasing property ${\cal Q}_n=\{\mbox{hypergraph is not $r$-colorable}\}$, where $r=r(n)\ge 2$ is some function of $n$. In the next paragraph we shall give a background of this problem.

\subsection{Previous results}

The $r$-colorability of random hypergraph $H(n,k,p)$ was most intensively studied in the case of fixed $k$ and $r=2$. It appears that in this case the transition from 2-colorability to non-2-colorability is \emph{sharp}. It follows form the general results of E. Freidgut (see \cite{Friedgut}) that for any $k\ge 3$, there exists a sequence $d_k(n)$ such that for any $\varepsilon>0$,
\begin{itemize}
  \item if $p\le (d_k(n)-\varepsilon)n/{n\choose k}$, then ${\sf P}\left(H(n,k,p)\mbox { is 2-colorable}\right)\to 1$,

  \item but if $p\ge (d_k(n)+\varepsilon)n/{n\choose k}$, then ${\sf P}\left(H(n,k,p)\mbox { is 2-colorable}\right)\to 0$.
\end{itemize}
It is widely believed that $d_k(n)$ can be replaced by a constant $d_k$.

First bounds for the threshold probability of 2-colorability were obtained by N. Alon and J. Spencer. They showed (see \cite{AlonSpencerUnp}) that there is a positive absolute constant $c$ such that
\begin{equation}\label{alon_spencer_01}
\mbox{ if }k\ge k_0\mbox{ is fixed }\mbox{ and }p\le c\,\frac{2^{k-1}}{k^2}\,\frac {n}{{n\choose k}},\mbox{ then }{\sf P}\left(H(n,k,p)\mbox { is 2-colorable}\right)\to 1,
\end{equation}
$$
\mbox{ if }k\ge 3,\;\varepsilon>0\mbox{ are fixed }\mbox{ and }p\ge (1+\varepsilon)\,2^{k-1}\ln 2\,\frac {n}{{n\choose k}},
$$
\begin{equation}\label{alon_spencer_02}
\mbox{ then }{\sf P}\left(H(n,k,p)\mbox { is 2-colorable}\right)\to 0.
\end{equation}

The gap between upper and lower bounds in (\ref{alon_spencer_01}) and (\ref{alon_spencer_02}) was reduced by D. Achlioptas, J.H. Kim, M. Krivelevich and P. Tetali from the order $k^2$ to the order $k$. They proved (see \cite{AKKT}) that for any fixed $k\ge 3$,
\begin{equation}\label{akkt_2colors}
  \mbox{ if }p\le \frac 1{25}\,\frac {2^{k-1}}{k}\,\frac {n}{{n\choose k}},\mbox{ then }{\sf P}\left(H(n,k,p)\mbox { is 2-colorable}\right)\to 1.
\end{equation}
Finally, Achlioptas and C. Moore established (see \cite{AchMoore}) the following bound for the threshold probability of 2-colorability for all sufficiently large $k$: there is a constant $k_0$ such that, for any $\varepsilon>0$ and any fixed $k\ge k_0$,
\begin{equation}\label{ach_moore}
\mbox{ if }p\le (1-\varepsilon)\,2^{k-1}\ln 2\,\frac {n}{{n\choose k}},\mbox{ then }{\sf P}\left(H(n,k,p)\mbox { is 2-colorable}\right)\to 1.
\end{equation}
Together with the upper bound (\ref{alon_spencer_02}) the inequality (\ref{ach_moore}) gives the exact value of the considered sharp threshold for 2-colorability in the case of fixed $k\ge k_0$:
$$
  p^*=2^{k-1}\ln 2\,\frac {n}{{n\choose k}}.
$$

\bigskip
The $r$-colorability of $H(n,k,p)$ for $r>2$ is not studied in such detail as 2-colorability. The following lemmas are just natural generalizations of the results (\ref{alon_spencer_01}) and (\ref{alon_spencer_02}) of Alon and Spencer.

\begin{lemma}\label{lemma_1} There exist positive constants $C,c>0$ such that for any $k=k(n)\ge 3$ and $r=r(n)\ge 2$, satisfying the conditions $r^{k-1}/k\ge C$ and $r^{k-1}=o(n)$ the following statement holds:
\begin{equation}\label{sta_lemma_1}
\mbox{ if }p\le c\,\frac{r^{k-1}}{k^2}\,\frac {n}{{n\choose k}},\mbox{ then }{\sf P}\left(H(n,k,p)\mbox { is $r$-colorable}\right)\to 1.
\end{equation}
\end{lemma}
\vspace{2mm}
\begin{lemma}\label{lemma_2} Let the functions $k=k(n)$ and $r=r(n)$ satisfy the relation $k^2r=o(n)$. Then for any fixed $\varepsilon>0$,
\begin{equation}\label{sta_lemma_2}
\mbox{ if }p\ge (1+\varepsilon)\,r^{k-1}\ln r\,\frac {n}{{n\choose k}},\mbox{ then }{\sf P}\left(H(n,k,p)\mbox { is $r$-colorable}\right)\to 0.
\end{equation}
\end{lemma}

Another result for $r$-colorability of random hypergraphs was obtained by Achlioptas, Kim, Krivelevich and Tetali. In the final comment of their paper \cite{AKKT} they stated (with providing an algorithm of the proof) that the result (\ref{akkt_2colors}) can be generalized to the case of $r$ colors in the following form.

\begin{theorem}\label{thm_akkt_rcolors}{\rm (D. Achlioptas, J.H. Kim, M. Krivelevich, P. Tetali, \cite{AKKT})} Suppose $k\ge 3$ and $r\ge 2$ are fixed. If $p=p(n)$ satisfies the inequality
\begin{equation}\label{akkt_rcolors}
  p\le \frac {r(r+1)!}{(r+1)^{2(r+1)}}\,\frac {r^{k-1}}{k}\,\frac {n}{{n\choose k}},
\end{equation}
then ${\sf P}\left(H(n,k,p)\mbox { is r-colorable}\right)\to 1$.
\end{theorem}
It is easy to see that the bound (\ref{sta_lemma_1}) of Lemma \ref{lemma_1} is better than (\ref{akkt_rcolors}) if
$$
  k\le c \frac {(r+1)^{2(r+1)}}{r(r+1)!}=\Omega\left(r^{r}\right).
$$
So, Theorem \ref{thm_akkt_rcolors} gives a new result only when $r$ is small in comparison with $k$: $r=O(\ln k/\ln\ln k)$.

\bigskip
The threshold probability for $r$-colorability of the random hypergraph $H(n,k,p)$ in the case when $r$ is large in comparison with $k$ can be obtained by using the results concerning the chromatic number of $H(n,k,p)$ (recall, e.g., that (\ref{ach_moore}), (\ref{akkt_rcolors}) are nontrivial only when $k$ is much larger than $r$). This problem was studied by a series of researchers (see, e.g., \cite{KarLuc}, \cite{KrivSud} for the background). In our paper we study $H(n,k,p)$ in the ``sparse'' case, i.e. the function $p=p(n)$ is sufficiently small. For such values of $p$, M. Krivelevich and B. Sudakov proved (see \cite{KrivSud}) the following theorem.

\begin{theorem}\label{thm_krivsud}{\rm (M. Krivelevich, B. Sudakov, \cite{KrivSud})} Let $k\ge 3$ be fixed. There is a constant $d_0=d_0(k)$ such that, for any $p=p(n)$ satisfying the conditions
$$
  d=d(n)=(k-1){n-1\choose k-1}p\ge d_0,\quad d=o(n^{k-1}),
$$
the following convergence holds:
$$
  {\sf P}\left(\left(\frac{d}{k\ln d}\right)^{1/(k-1)}\le \chi(H(n,k,p))\le \left(\frac{d}{k\ln d}\left(1+\frac {28k\ln\ln d}{\ln d}\right)\right)^{1/(k-1)}\right)\to 1.
$$
\end{theorem}

One can make an immediate corollary from this theorem.

\begin{corollary}\label{cor_krivsud} Let $k\ge 3$ and $\varepsilon\in(0,1)$ be fixed. There is a constant $r_0=r_0(k,\varepsilon)$ such that, for any $r=r(n)$ satisfying the conditions
$$
  r\ge r_0,\quad r^{k-1}\ln r=o(n^{k-1}),
$$
the following convergence holds:
$$
  {\sf P}\left(\chi(H(n,k,p))\le r\right)\to 1,\mbox{ where }p=(1-\varepsilon)\,r^{k-1}\ln r\,\frac{n}{{n\choose k}}.
$$
\end{corollary}
Corollary \ref{cor_krivsud} together with Lemma \ref{lemma_2} shows that the function $p^*=r^{k-1}\ln r\,n/{n\choose k}$ is a threshold probability in the wide range: $k$ is fixed, $r$ is sufficiently large in comparison with $k$ and $r^{k-1}\ln r=o(n^{k-1})$.

However, Theorem \ref{thm_krivsud} (and, consequently, Corollary \ref{cor_krivsud}) can be proved not only for fixed $k$, but for slowly growing functions $k=k(n)$ also. The calculations from the proof of Theorem \ref{thm_krivsud} provides the following necessary relations between $d$, $p$, $k$ and $n$:
\begin{equation}\label{ks_condition}
  d\ge (\ln d)^{28k-27},\quad n^{1/3}\ge\left(\ln (n^{k-1}p)\right)^{3(k-1)-1/2}.
\end{equation}
These relations implies that in Corollary \ref{cor_krivsud} we have the following restrictions:
\begin{equation}\label{cor_ks_condition}
  r=\Omega\left(k^{29}(\ln k)^{28}\right),\quad n\ge k^{9k+O(k\ln\ln k/(\ln k))}.
\end{equation}
So, despite the fact that Corollary \ref{cor_krivsud} gives very good lower bound for the threshold probability, its statement holds only for large $r$ in comparison with $k$: $r=\Omega\left(k^{29}(\ln k)^{28}\right)$. Recall that (\ref{akkt_rcolors}) is better than (\ref{sta_lemma_1}) only when $r=O(\ln k/\ln\ln k)$. Hence, in the very wide range of the values of $r$,
\begin{equation}\label{area}
  \frac {\ln k}{\ln\ln k}\le r\le k^{29}(\ln k)^{28},
\end{equation}
only the lower bound from Lemma \ref{lemma_1} is known.

\bigskip
\begin{note}
The proof of Theorem \ref{thm_krivsud} seems possible to be adopted to the case of smaller values of the parameter $d$ (and, consequently, parameter $r$ in Corollary \ref{cor_krivsud}) than given by (\ref{ks_condition}) . But, for example, the final condition $r>k^4$ seems to be necessary. So, the case when $r$ is not very large in comparison with $k$ is certainly not well studied.
\end{note}

\bigskip
We have finished discussing previously known results and now proceed to the new ones.

\section{New results}

Our main approach of studying the threshold for $r$-colorability of random hypergraph $H(n,k,p)$ is to apply methods and results concerning extremal problems of hypergraph coloring theory.

\subsection{Colorings of hypergraphs with bounded vertex degrees}

For all $k,r\ge 2$, let $\Delta(k,r)$ denote the minimum possible $\Delta(H)$, where $H$ is a $k$-uniform non-$r$-colorable hypergraph. The problem of finding or estimating the value $\Delta(k,r)$ is one of the classical problems in extremal combinatorics. First bounds for $\Delta(k,r)$ were obtained by P. Erd\H{o}s and L. Lov\'asz (see \cite{ErdLov}), they proved that for all $k,r\ge 2$,
\begin{equation}\label{d(k,r)_erdlov}
  \frac {r^{k-1}}{4k}\le \Delta(k,r)\le 20k^2r^{k+1}.
\end{equation}
Kostochka and R\"odl improved (see \cite{KostRodl}) the upper bound from (\ref{d(k,r)_erdlov}). They showed that for all $k,r\ge 2$,
$$
	\Delta(k,r)\le \left\lceil k\,r^{k-1}\ln r\right\rceil.
$$

Classical lower bound (\ref{d(k,r)_erdlov}) of Erd\H{o}s and Lov\'asz was improved by J. Radhakrishnan and A. Srinivasan (see \cite{RadhSrin}) in the case $r=2$. They proved that for large $n$,
$$
	\Delta(k,2)\ge 0.17\, \frac {2^{k}}{\sqrt{k\ln k}}.
$$
Their result is still the best one in the case of two colors.

When $r>2$ D.A. Shabanov proved (see \cite{Shab1}) a lower bound with slightly better ``polynomial'' factor: for any $k\ge 3$, $r\ge 3$,
\begin{equation}\label{d(k,r)_shab}
	\Delta(k,r)>\frac 18\, k^{-1/2}r^{k-1}.
\end{equation}
The last known result concerning $\Delta(k,r)$ was recently obtained by A.V. Kostochka, M. Kumbhat and V. R\"odl (see \cite{KKR}). They showed that if $r=r(n)\ll\sqrt{\ln\ln k}$, then
\begin{equation}\label{d(k,r)_kkr}
\Delta(k,r)>e^{-4r^2} \left(\frac k{\ln k}\right)^{\frac {\left\lfloor \log_2 r\right\rfloor}{\left\lfloor \log_2 r\right\rfloor+1}}\frac{r^k}k.
\end{equation}

\bigskip
The following lemma clarifies the connection between the value $\Delta(k,r)$ and the threshold for $r$-colorability of random hypergraph $H(n,k,p)$.
\begin{lemma}\label{lemma_3} Suppose $k=k(n)\ge 2$ and $r=r(n)\ge 2$ satisfy the relation
\begin{equation}\label{cond_lemma3}
  \frac 3{16}\,\Delta(k,r)-\ln n\to -\infty\mbox{ as }n\to\infty.
\end{equation}
If
$$
  p\le \frac 12\,\frac{\Delta(k,r)}{k}\,\frac n{{n\choose k}},
$$
then ${\sf P}\left(H(n,k,p)\mbox { is $r$-colorable}\right)\to 1$.
\end{lemma}

\begin{proof1} Since the probability ${\sf P}\left(H(n,k,p)\mbox { is $r$-colorable}\right)$ decreases with growth of $p$, we have to deal only with $p=\frac 12\,\frac{\Delta(k,r)}{k}\,\frac n{{n\choose k}}=\frac 12\,\Delta(k,r){n-1\choose k-1}^{-1}$.

Let $v$ be a vertex of $H(n,k,p)$ and let $X_v$ denote the degree of $v$ in $H(n,k,p)$. It is clear that $X_v$ is a binomial random variable $Bin\left({n-1\choose k-1},p\right)$. We shall need a classical bound on probability of large deviations for binomial variables (so called, Chernoff bound): if $X$ is a binomial random variable, then for any $\lambda>0$,
\begin{equation}\label{chernoff}
	{\sf P}\left(X\ge {\sf E}X+\lambda\right)\le \exp\left\{-\frac{\lambda^2}{2({\sf E}X+\lambda/3)}\right\}.
\end{equation}
The proof of this classical fact can be found, e.g., in the book \cite{JLR}. Using (\ref{chernoff}) with $\lambda={\sf E}X$ we get
$$
	{\sf P}\left(X_v\ge 2 {n-1\choose k-1}p\right)={\sf P}\left(X_v\ge \Delta(k,r)\right)\le
  \exp\left\{-\,\frac 3{16}\,\Delta(k,r)\right\}.
$$

Consequently we obtain the following bound for the probability of the existence of the vertices with large degree in $H(n,k,p)$:
$$
  {\sf P}\left(\Delta(H(n,k,p))\ge \Delta(k,r)\right)\le n\,\exp\left\{-\,\frac 3{16}\,\Delta(k,r)\right\}=\exp\left\{\ln n-\,\frac 3{16}\,\Delta(k,r)\right\}\to 0
$$
as $n\to+\infty$. The last relation follows from the condition (\ref{cond_lemma3}). Thus, by the definition of the value $\Delta(k,r)$ we have
$$
  \lim_{n\to\infty}{\sf P}(\chi(H(n,k,p))\le r)\ge \lim_{n\to\infty}{\sf P}\left(\Delta(H(n,k,p))<\Delta(k,r)\right)=1.
$$
Lemma \ref{lemma_3} is proved.
\end{proof1}

\bigskip
As a corollary of Lemma \ref{lemma_3} and the bounds (\ref{d(k,r)_shab}) and (\ref{d(k,r)_kkr}) for $\Delta(k,r)$ we immediately obtain the following lower bound for the threshold for $r$-colorability of $H(n,k,p)$.

\begin{corollary}\label{corollary1} 1) Suppose $k=k(n)\ge 3$ and $r=r(n)\ge 3$ satisfy the relation
\begin{equation}\label{cond1_corollary1}
  \frac{3}{128}\,\frac{r^{k-1}}{\sqrt{k}}-\ln n\to -\infty\mbox{ as }n\to\infty.
\end{equation}
If
\begin{equation}\label{thresh_bound1}
  p\le \frac 3{32}\,\frac{r^{k-1}}{k^{3/2}}\,\frac n{{n\choose k}},
\end{equation}
then ${\sf P}\left(H(n,k,p)\mbox { is $r$-colorable}\right)\to 1$.

2) Suppose $k=k(n)\ge 3$ and $r=r(n)\ge 2$ satisfy the relation
$$
  \frac{3}{16}\,e^{-4r^2}\left(\frac k{\ln k}\right)^{\frac {\left\lfloor \log_2 r\right\rfloor}{\left\lfloor \log_2 r\right\rfloor+1}}\frac{r^k}{k}-\ln n\to -\infty\mbox{ as }n\to\infty.
$$
If $r=o(\sqrt{\ln\ln k})$ and
\begin{equation}\label{thresh_bound2}
  p\le \frac 3{16}\,e^{-4r^2}\left(\frac k{\ln k}\right)^{\frac {\left\lfloor \log_2 r\right\rfloor}{\left\lfloor \log_2 r\right\rfloor+1}}\frac{r^k}{k^2}\,\frac n{{n\choose k}},
\end{equation}
then ${\sf P}\left(H(n,k,p)\mbox { is $r$-colorable}\right)\to 1$.
\end{corollary}

Let us compare the results of Corollary \ref{corollary1} with previous ones. Since Corollary \ref{cor_krivsud} gives almost complete answer, we have to compare (\ref{thresh_bound1}) and (\ref{thresh_bound2}) with (\ref{sta_lemma_1}) from Lemma \ref{lemma_1} and (\ref{akkt_rcolors}) from Theorem \ref{thm_akkt_rcolors}. Although the bounds (\ref{sta_lemma_1}) and (\ref{akkt_rcolors}) formally hold only when $k$ and $r$ are fixed, the analogous statements can be proved by using the same arguments for growing functions $k=k(n)$ and $r=r(n)$. For example, the analysis of the calculations in the papers \cite{AlonSpencerUnp}, \cite{AKKT} and \cite{AchMoore} shows that
\begin{enumerate}
  \item The statement of Lemma \ref{lemma_1} holds for almost all functions $k=k(n)$ and $r=r(n)$, since it is true for $r^{k-1}=o(n)$ and in the case $r^{k-1}/k\ge 22\ln n$ we can just apply Lemma \ref{lemma_3} with the classical bound (\ref{d(k,r)_erdlov}) of Erd\H{o}s and Lov\'asz.
  \item The proof of Theorem \ref{thm_akkt_rcolors} can be extended to the following range of values of $k=k(n)$ and $r=r(n)$: $r$ is fixed and $k=o(\sqrt{n})$.
  \item The proof of the result (\ref{ach_moore}) by Achlioptas and Moore does not work for any growing $k=k(n)$. It is also unclear how to generalize it to case of fixed $r>2$. Thus, we do not compare our new results with (\ref{ach_moore}), since we consider only the situation when $k$ or $r$ (or both of them) is a growing function of $n$.
\end{enumerate}
Let us sum up the obtained information. The lower bound (\ref{sta_lemma_1}) from Lemma \ref{lemma_1} holds almost for all $r$ and $k$. Theorem \ref{thm_akkt_rcolors} also can be extended to a wide area of the values of the parameters. Everywhere below for simplicity we shall compare only the values of the bounds.

Both (\ref{thresh_bound1}) and (\ref{thresh_bound2}) are obviously better than (\ref{sta_lemma_1}). The second bound (\ref{thresh_bound2}) is worse than (\ref{akkt_rcolors}). Indeed, the right hand-side of (\ref{akkt_rcolors}) is at least
$$
  e^{-2r\ln r}\,\frac{r^{k-1}}{k},
$$
which is better than (\ref{thresh_bound2}), whose right hand-side is at most
$$
  e^{-4r^2}\frac{r^{k-1}}{k^{(\lfloor\log_2 r\rfloor+2)/(\lfloor\log_2 r\rfloor+1)}}.
$$
The first bound (\ref{thresh_bound1}) of Corollary \ref{corollary1} is better than (\ref{akkt_rcolors}) if
$$
  \sqrt{k}<\frac 3{32}\,\frac{(r+1)^{2(r+1)}}{r(r+1)!}.
$$
This inequality holds, e.g., when $r\ge \ln k/\ln\ln k$.

\bigskip
Let us make intermediate conclusions. Our new lower bound (\ref{akkt_rcolors}) for the threshold probability of $r$-colorability of random hypergraph $H(n,k,p)$ improves all previously known results in the following wide area (see condition (\ref{cor_ks_condition}) of Corollary \ref{cor_krivsud} and condition (\ref{cond1_corollary1}) of Corollary \ref{corollary1}):
\begin{equation}\label{bound1_area}
  \ln k/\ln\ln k\le r\le k^{29}(\ln k)^{28}\quad\mbox { and }\quad\frac{r^{k-1}}{\sqrt{k}}\gg \ln n.
\end{equation}
We see that in the area (\ref{bound1_area}) the parameter $r$ cannot be very small in comparison with the number of vertices $n$, but it can be very large. In the next paragraph we shall present a better bound when $r$ is not very small and also is not very large in comparison with $n$.

\subsection{Main result}

The main result of our paper is formulated in the following theorem.

\begin{theorem}\label{theorem 2}
Suppose $\delta\in(0,1)$ is a constant. Let $k=k(n)$ and $r=r(n)\ge 2$ satisfy the following conditions: $k\ge k_0(\delta)$ , where $k_0(\delta)$ is some constant, and, moreover,
\begin{equation}\label{n th 2}
(k-1)\ln r < \frac {1-\delta}2\ln n,\quad r^{k-1}k^{-\varphi(k)}\ge 6\ln n,
\end{equation}
where $\varphi(k)=4\left\lfloor\sqrt{\frac{\ln k}{\ln(2\ln k)}}\right\rfloor^{-1}$. Then for  function $p=p(n)$, satisfying
\begin{equation}\label{p th 2}
p\le \frac 12\,\frac {r^{k-1}}{k^{1+\varphi(k)}}\,\frac n{{n\choose k}},
\end{equation}
we have ${\sf P}\left(H(n,k,p)\mbox{ is r-colorable}\right)\to 1$ as $n\to\infty$.
\end{theorem}

Let us compare the result of Theorem \ref{theorem 2} with previous ones. It is clear that the restriction (\ref{p th 2}) is weaker (for all sufficiently large $k$) than our previous results (\ref{thresh_bound1}) and (\ref{thresh_bound2}) obtained in $\S$2.1. It is also obvious that (\ref{p th 2}) is better than the lower bound (\ref{sta_lemma_1}) from Lemma \ref{lemma_1} for all sufficiently large $k$. So, it remains only to compare (\ref{p th 2}) with the result of Theorem \ref{thm_akkt_rcolors} proved by Achlioptas, Krivelevich, Kim and Tetali. We see that (\ref{p th 2}) is asymptotically better than (\ref{akkt_rcolors}) if
$$
  k^{4\left\lfloor\sqrt{\frac{\ln k}{\ln(2\ln k)}}\right\rfloor^{-1}}<\frac{(r+1)^{2(r+1)}}{2r(r+1)!}.
$$
This inequality holds, e.g., in the following asymptotic area: $r\gg \sqrt{\ln k}$.

Thus, our main result (\ref{p th 2}) gives a new lower bound for the threshold probability for $r$-colorability of the random hypergraph $H(n,k,p)$. This new bound improves all the previously known results in the wide area of the parameters (recall that we are working in the area (\ref{area})):
$$
  \sqrt{\ln k}\ll r\le k^{29}(\ln k)^{28}\;\mbox{ and }\; 6\,k^{\varphi(k)}\ln n\le r^{k-1}\le n^{(1-\delta)/2}.
$$
For example, (\ref{p th 2}) provides a new bound when $k\sim r\sim \ln n/(5\ln\ln n)$. Moreover, our result (\ref{p th 2}) is only $k^{1+\varphi(k)}\ln r$ times smaller than the upper bound (\ref{sta_lemma_2}) from Lemma \ref{lemma_2}.

\bigskip
The proof of Theorem \ref{theorem 2} is based on some result concerning colorings of 2-simple hypergraphs with bounded vertex degrees. The study of problems for colorings of simple hypergraphs was initiated by Erd\H{o}s and Lov\'asz in \cite{ErdLov}. Later the extremal problems concerning colorings of $l$-simple hypergraphs with bounded vertex degrees have been considered by Z. Szab\'o (see \cite{Szabo}), A.V. Kostochka and M. Kumbhat (see \cite{KostKumb}), D.A. Shabanov (see \cite{Shab2}). To prove Theorem \ref{theorem 2} we consider 2-simple hypergraphs with a few 3-cycles.

Let $H$ be an arbitrary hypergraph with the following properties: $H$ is $k$-uniform, $\chi(H)>r$, $H$ is 2-simple and for every edge of $H$, there are at most $\omega$ 3-cycles containing that edge. The class of all such hypergraphs we will denote by $\mathcal{H}(k,r,\omega).$ Let us consider the following extremal value:
$$
  \Delta(\mathcal{H}(k,r,\omega))=\min\left\{\Delta(H):\;H\in \mathcal{H}(k,r,\omega)\right\}.
$$
Theorem \ref{theorem 3} gives an asymptotic lower bound for $\Delta(\mathcal{H}(k,r,\omega))$ for $\omega=\lfloor \sqrt{\ln k/(\ln\ln k)}\rfloor$.
\begin{theorem}\label{theorem 3}
There exists an integer $k_0$ such that for all $k\ge k_0$, all $r\ge 2$ and all $\omega\le \sqrt{\ln k/(\ln\ln k)}$,
\begin{equation}\label{delta_low}
	\Delta(\mathcal{H}(k,r,\omega))>r^{k-1}k^{-4\left\lfloor\sqrt{\frac{\ln k}{\ln(2\ln k)}}\right\rfloor^{-1}}.
\end{equation}
\end{theorem}

It should be noted that the inequality (\ref{delta_low}) holds for all possible values of the parameter $r$, which is important for studying $r$-colorability of random hypergraphs. For a $k$-uniform, 2-simple, non-$r$-colorable hypergraph the lower bound for the maximum vertex degree similar to (\ref{delta_low}) is known only in the case of small $r$ in comparison with $k$: $r=O(\ln k)$ (see \cite{Shab2} for the details).

\bigskip
The structure of the rest of the article will be the following. In the next paragraph we shall deduce Theorem \ref{theorem 2} from Theorem \ref{theorem 3}. Section 3 will be devoted to the proof of Theorem \ref{theorem 3}. Finally, in Section 4 we shall discuss choosability in random hypergraphs.

\subsection{Proof of Theorem \ref{theorem 2}}

Due to the decreasing of the probability ${\sf P}\left(H(n,k,p)\mbox { is $r$-colorable}\right)$ with growth of $p$ we have to deal only with $p=\frac {r^{k-1}}{2k^{1+\varphi(k)}}\,\frac n{{n\choose k}}$. We want to apply Theorem \ref{theorem 3} to random hypergraph $H(n,k,p)$, so, we have to show that with probability tending to 1 $H(n,k,p)$ satisfies the following conditions: it is 2-simple, every edge is contained in at most $\omega=\lfloor\sqrt{\ln k/\ln\ln k}\rfloor$ 3-cycles and, moreover, $\Delta(H(n,k,p))<\Delta(\mathcal{H}(k,r,\omega))$.

Let $v$ be a vertex of the random hypergraph $H(n,k,p)$ and let $X_v$ denote the degree of $v$ in $H(n,k,p)$. It is clear that $X_v$ is a binomial random variable $Bin\left({n-1\choose k-1},p\right)$. Using Chernoff bound (\ref{chernoff}) with $\lambda={\sf E}X_v$ and the condition (\ref{n th 2}), we have
$$
  {\sf P}\left(X_v\ge r^{k-1}k^{-\varphi(k)}\right)\le \exp\left\{-3r^{k-1}k^{-\varphi(k)}/16\right\}\le\exp\{-(9/8)\ln n\}=n^{-9/8}.
$$
Thus,
\begin{equation}\label{thm2_ineq1}
  {\sf P}\left(\Delta(H(n,k,p))\ge r^{k-1}k^{-\varphi(k)}\right)\le n\cdot n^{-9/8}=o(1).
\end{equation}

Let $Y$ denote the number of pairs of edges, whose intersection has cardinality at least 3, and  let $Z$ denote the number of edges, which are contained in a large number, more than $\omega$, of 3-cycles. We estimate the expectations of these two random variables:
$$
  {\sf E}Y \le {n\choose 3}{n-3\choose k-3}^2 p^2\le \frac {n^3k^6}{n^6}{n\choose k}^2\left(\frac {r^{k-1}}{k}n {n\choose k}^{-1}\right)^2=\frac {\left(r^{k-1}\right)^2k^4}n.
$$
Consequently,
$$
  \ln {\sf E}Y = (2k-2)\ln r+4\ln k-\ln n\overset{(\ref{n th 2})}{\le} (1-\delta)\ln n +4\ln k-\ln n\le -\frac{\delta}2 \ln n.
$$
Hence, $\lim_{n\to \infty} {\sf E}Y=0$ and
\begin{equation}\label{thm2_ineq2}
  {\sf P}\left(H(n,k,p)\mbox{ is 2-simple}\right)\to 1.
\end{equation}

Now we will consider edges, that are contained in a large number of triangles. Suppose $u$ is an edge of $H(n,k,p)$. Let us denote by $T_u$ the set of all triangles, containing $u$. Furthermore,
we denote by $D(u',u)$ \textit{the degree of an edge $u'$ with respect to $u$}, a number of 3-cycles from $T_u$, containing an edge $u'\neq u$. Similarly, for any vertex $v\in V(H(n,k,p))$, we denote by $d(v,u)$ the \textit{the degree of vertex $v$ with respect to $u$}, a number of triangles $(u,u',u'')$ from $T_u$ such that $v\in (u'\cap u'')\setminus u$.

Now we will estimate the number of edges that have big degree with respect to $T_u$ for some $u$.
Denote by $Z_1(u)$ the number of edges, having degree greater than 4 with respect to $u$, i.e.
$$
  Z_1(u)=\left|\left\{u'\in E(H(n,k,p)):\;D(u',u)>4\right\}\right|.
$$
Moreover, let us denote $Z_1=\sum\limits_{e\in\;E(H(n,k,p)}Z_1(u)$. Now we estimate the expectation of $Z_1$:
$$
   {\sf E}Z_1\le n {n-1\choose k-1}^2\,k^{10}{n-2\choose k-2}^5p^7\le
   {n\choose k}^7\,\frac {k^{22}}{n^{11}}\left(\frac {r^{k-1}}{k}n {n\choose k}^{-1}\right)^7=\frac{r^{7(k-1)}k^{15}}{n^{4}}=
$$
$$
  =\exp\left\{7(k-1)\ln r+15\ln k-4\ln n\right\}\overset{(\ref{n th 2})}{\le}
  \exp\left\{\frac 72(1-\delta)\ln n +15\ln k-4\ln n\right\}\le n^{-1/2}.
$$
Let us explain the first inequality. At first we choose the vertex from the intersection of the edge $u'$ with large degree and the edge $u$. Then we choose the rest vertices of these two edges. Then we choose 5 vertices on both of edges, that correspond to remaining vertices of five 3-cycles. Then we choose the last edge of each 3-cycle.

Thus,
\begin{equation}\label{thm2_ineq3}
  {\sf P}\left(\mbox{for any }u,u'\in E(H(n,k,p)),\;\;D(u',u)\le 4\right)\to 1.
\end{equation}

Similarly, we shall show, that with probability tending to one, $d(v,u)\le 4$ for any edge $u$ and any vertex $v\notin u$. Namely, we denote by $Z_2$ the number of pairs $v,u$, such that $d(v,u)\ge 5$. Then the expectation of $Z_2$ can be estimated from above as follows.
$$
  {\sf E}Z_2\le n {n-1\choose k}\,k^{5}{n-2\choose k-2}^5p^6\le
  {n\choose k}^6\,\frac {k^{15}}{n^{9}}\left(\frac {r^{k-1}}{k}n {n\choose k}^{-1}\right)^6=\frac{r^{6(k-1)}k^{9}}{n^{3}}=
$$
$$
=\exp\left\{6(k-1)\ln r+9\ln k-3\ln n\right\}\overset{(\ref{n th 2})}{\le}
\exp\left\{3(1-\delta)\ln n +9\ln k-3\ln n\right\}\le n^{-\delta/2}.
$$
So,
\begin{equation}\label{thm2_ineq4}
  {\sf P}\left(\mbox{for any vertex }v\mbox{ and an edge }u\in E(H(n,k,p)),\;\;d(v,u)\le 4\right)\to 1.
\end{equation}

Let us introduce the following event
$$
  {\cal A}_n=\left\{\mbox{for any }u,u'\in E(H(n,k,p))\mbox{ and any }v\in V(H(n,k,p)),\;D(u',u)\le 4\mbox{ and }d(v,u)\le 4\right\}.
$$
Due to (\ref{thm2_ineq3}) and (\ref{thm2_ineq4}) we have that ${\sf P}({\cal A}_n)\to 1$.

Now suppose that the event ${\cal A}_n$ holds and there is an edge $u$ in $H(n,k,p)$, which is contained in at least $\omega$ 3-cycles. Consider the following set of vertices $V_u$:
$$
  V_u=\{v\in V:\;\; d(v,u)>0\}.
$$
For any $v\in V_u$, by $E(v,u)$ we denote the set of edges, containing $v$, which also belongs to one of the 3-cycles from $T_u$. The event ${\cal A}_n$ implies that $d(v,u)\le 4$ for any $v$ and $u$ and, hence, $|E(v,u)|\le 8$ and $|V_u|\ge \frac 18|T_u|=\frac 18\omega$. Now we will construct a some convenient subset  of $T_u$ of sufficient size.

First, we have $V_u^0=V_u$ and $T_u^0=T_u$. Suppose sets $V_u^s$ and $T_u^s$ are considered. We form a set $V_u^{s+1}$ and a set $T^{s+1}_u$ by the following way. We choose an arbitrary 3-cycle $t_{s+1}=(u^{'}_s,u^{''}_s,u)\in T_u^s$ and an arbitrary $v_{s+1}\in (u^{'}_s\cap u^{''}_s)\setminus u$. Then we take
$$
  V_u^{s+1}=V_u^s\setminus\{v_{s+1}\},
$$
$$
  T_u^{s+1}=T_u^s\setminus\{t\in T_u^s:\;t\mbox{ contains an edge from }E(v_{s+1},u)\}.
$$
Then we repeat the same procedure with sets $T_u^{s+1}, V_u^{s+1}$. The procedure continues until both sets $V_u^{s+1}$ and $T_u^{s+1}$ are not empty.

How many steps of procedure can we guarantee? Since the event ${\cal A}_n$ holds, we have that$|T^{i+1}_u|\ge |T^{i}_u|-32$. Indeed, $E(v_{i+1},u)$ consists of at most 8 edges and every edge belongs to at most 4 3-cycles. So, we can guarantee at least $\omega'=\lceil\frac 1 {32}\omega\rceil$ steps of the procedure.

Consider the obtained set of 3-cycles $\{t_1,\ldots,t_{\omega'}\}$, $t_s=(u^{'}_s,u^{''}_s,u)$, $s=1,\ldots,\omega'$ and the set of vertices $\{v_1,\ldots,v_{\omega'}\}$. Our procedure shows that all the edges $u^{'}_1,u^{''}_1,u^{'}_2,u^{''}_2,\ldots,u^{'}_{\omega'},u^{''}_{\omega'}$ are distinct, all the vertices $v_1,\ldots,v_{\omega'}$ are also distinct and, for any $s=1,\ldots,\omega'$, we have $v_s\in u^{'}_s\cap u^{''}_s\setminus u$. Let us estimate the probability of the event (denoted by ${\cal B}_n$) that, for some edge $u$, the described above configuration of appears in $H(n,k,p)$. It is clear that
$$
  {\sf P}\left({\cal B}_n\right)\le {n\choose k}k^{2\omega'}(n-k)^{\omega'}{n-2\choose k-2}^{2\omega'} p^{2\omega'+1}
  \le \frac {n^{\omega'}k^{6\omega'}}{n^{4\omega'}}{n\choose k}^{2\omega'+1}\left(\frac {r^{k-1}}{k}n {n\choose k}^{-1}\right)^{2\omega'+1}=
$$
$$
  =\frac {\left(r^{k-1}\right)^{2\omega'+1}k^{4\omega'-1}}{n^{\omega' -1}}.
$$
Hence, for all $k\ge k_0$,
$$
  \ln {\sf P}\left({\cal B}_n\right) = (k-1)(2\omega'+1)\ln r+(4\omega'-1)\ln k-(\omega'-1)\ln n\overset{(\ref{n th 2})}{\le}
$$
\begin{equation}\label{thm2_ineq6}
  \overset{(\ref{n th 2})}{\le} (\omega'-1)\left((1-\delta)\left(\frac {2\omega'+1}{2\omega'-2}\right)\ln n +\frac{4\omega'-1}{(\omega'-1)}\ln k-\ln n\right)\le(\omega'-1)\left( -\frac{\delta}2 \ln n\right).
\end{equation}
Consequently, $\lim_{n\to \infty} {\sf P}\left({\cal B}_n\right)=0$. Finally, if ${\cal A}_n$ holds, then the event that there is an edge $u\in E(H(n,k,p))$ with $|T_u|>\omega$ implies the evnt ${\cal B}_n$. Thus,
\begin{equation}\label{thm2_ineq7}
  {\sf P}\left(\mbox{there is }u\in E(H(n,k,p)):\;|T_u|>\omega\right)\le {\sf P}\left(\overline{{\cal A}_n}\right)+{\sf P}\left({\cal B}_n\right)\to 0.
\end{equation}

\bigskip
Let us sum up the obtained information. It follows from (\ref{thm2_ineq2}) and (\ref{thm2_ineq7}) that with probability tending to 1 the random hypergraph $H(n,k,p)$ satisfies all the conditions of being an element of $\mathcal{H}(k,r,\omega(k)),$ except the condition $\chi(H(n,k,p))>r$. Applying Theorem \ref{theorem 3} and (\ref{thm2_ineq1}) we get that $H(n,k,p)\notin \mathcal{H}(k,r,\omega(k))$ with high probability. Thus, ${\sf P}\left(H(n,k,p)\mbox{ is $r$-colorable}\right)\to 1$ as $n\to\infty$. Theorem \ref{theorem 2} is proved.

\vspace{4mm}
\begin{note} If $k=k(n)\to +\infty$ then the parameter $\delta\in(0,1)$ from Theorem \ref{theorem 2} can be taken equal to some infinitesimal function. For example, it follows from  (\ref{thm2_ineq6}) that
$$
  \delta= 50\sqrt{\frac{\ln\ln k}{\ln k}}
$$
is sufficient.
\end{note}

\section{Proof of Theorem \ref{theorem 3}}

The proof of Theorem \ref{theorem 3} is based on the method of random recoloring. This method in the case of two colors was developed in the papers of J. Beck \cite{Beck}, J. Spencer \cite{Spencer}, Radhakrishnan and Srinivasan \cite{RadhSrin}. In this paper we follow the work \cite{Shab2} concerning $r$-colorings of $l$-simple hypergraphs with bounded edge degrees.

The structure of the proof will be the following. In the next section we will formulate a multiparametric Theorem \ref{theorem 4} which provides a new lower bound of the maximum edge degree in a hypergraph from the class ${\cal H}(k,r,\omega)$. Then we will prove Theorem \ref{theorem 4}. Finally we deduce Theorem \ref{theorem 3} from Theorem \ref{theorem 4} by choosing the values of the required parameters.

\subsection{General theorem}

Theorem \ref{theorem 3} is a simple corollary of the following multiparametric theorem.

\begin{theorem}\label{theorem 4}
Let $k\ge 3$, $r\ge 2$, $\omega\ge 1$ be integers, let $b$, $\alpha$ be positive numbers. Let us denote:
\begin{equation}\label{t}
t = \left\lfloor \sqrt{\frac{\ln k}{\ln(\alpha\ln k)}}\right\rfloor,\quad
q = \frac{\alpha\ln k}k.
\end{equation}

Let $H=(V,E)$ be an $k$-uniform 2-simple hypergraph such that for every edge in $H$ there are at most $\omega$ 3-cycles that contain that edge. Let, moreover, every edge of hypergraph $H$ intersects at most $d$ other edges of $H$, where
\begin{equation}\label{d_def}
d\le r^{k-1}k^{1-b/t}-1.
\end{equation}
If the following inequalities hold
\begin{equation}\label{condition1}
b\le t < k-\omega,
\end{equation}
\begin{equation}\label{condition2}
\frac 2k\le q \le \frac 12,
\end{equation}
$$
 \frac {k^2}{2^k}+(t+1)k^{1-\alpha} e^{\alpha(\ln k) (t+\omega)/k}  \left(\alpha\ln k\right)^{t+\omega}+\frac {(t+1)^2}{t!}\, k^{2-b}(\alpha\ln k)^{t\omega}+
$$
\begin{equation}\label{condition3}
+(t+1)t\left(\frac {2e\alpha\ln k}{t-1}\right)^{t-1} k^{1+\alpha-b}<\frac 14
\end{equation}
then $\chi(H)\le r$.
\end{theorem}

The proof of this theorem is based on a method of vertex random coloring. To prove Theorem \ref{theorem 4} we have to show the existence of a proper vertex $r$-coloring for hypergraph $H$. We shall construct some random $r$-coloring and estimate the probability that this coloring is not proper for $H$. If this probability is greater than 0, then we prove the existence of a required coloring, and the theorem follows.

\subsection{Algorithm of random recoloring}

We follow the ideas of Radhakrishnan and Srinivasan from \cite{RadhSrin} and the construction from \cite{Shab2} concerning random recoloring. Let $V=\{v_1,\ldots,v_w\}$. The algorithm consists of two phases.

	{\bf Phase 1.} We color all vertices randomly and uniformly with $r$ colors, independently from each other. Let us denote the generated random coloring by $\chi_0$.
	
	The obtained coloring $\chi_0$ can contain monochromatic edges and ``almost monochromatic'' edges. An edge $e\in E$ is said to be \textit{almost monochromatic} in $\chi_0$ if there is a color $u$ such that
$$
	n-t-\omega+2\le |\{v\in e:\;\mbox{$v$ is colored by $u$ in $\chi_0$} \}|<n.
$$
In this case, the color $u$ is called \textit{dominating} in $e$. For every $v\in V$, $u=1,\ldots,r$, let us use the notations
$$
	{\cal M}(v)=\left\{e\in E:\; v\in e,\;\mbox{$e$ is monochromatic in $\chi_0$}\right\},
$$
$$
	{\cal AM}(v,u)=\left\{e\in E:\; v\in e,\;\mbox{$e$ is almost monochromatic in $\chi_0$ with dominating color $u$}\right\}.
$$

 {\bf Phase 2.} In this phase, we want to recolor some vertices from the edges, which are monochromatic in $\chi_0$.
We consider the vertices according to an arbitrary fixed order $v_1,\ldots,v_w$.
Let $\{\eta_1,\ldots,\eta_w\}$ be mutually independent equally distributed random variables, taking values $1,\ldots,r$
with the same probability $p$ (the value of the parameter $p$ will be chosen later) and the value $0$ with probability $1-rp$.
The recoloring procedure consists of $w$ steps.

\begin{enumerate}
	\item[Step 1.] Assume that ${\cal M}(v_1)\ne\emptyset$ and, moreover, there is no $u=1,\ldots,r$ and $e\in {\cal AM}(v_1,u)$ such that
	\begin{enumerate}
	\item $\eta_1=u$,
	\item $v_1$ is the only vertex in $e$, which is not colored by $u$ in $\chi_0$.
\end{enumerate}
	Then we try to recolor $v_1$ according to the value of the random variable $\eta_1$:
	\begin{itemize}
	\item if $\eta_1=0$, then we do not recolor $v_1$,
	\item if $\eta_1\ne 0$, then we recolor $v_1$ in the color $\eta_1$.
\end{itemize}
	In all the other situations, we do not change the color of $v_1$. Let $\chi_1$ be the coloring after considering $v_1$.
\end{enumerate}

Now let the vertices $v_1,\ldots,v_{i-1}$ have been considered, so that the coloring $\chi_{i-1}$ is obtained.

\begin{enumerate}
	\item[Step i.] Assume that some $f\in{\cal M}(v_i)$ is still monochromatic in $\chi_{i-1}$ and, moreover, there is no $u=1,\ldots,r$ and $e\in {\cal AM}(v_i,u)$ such that
	\begin{enumerate}
	\item $\eta_i=u$,
	\item $v_i$ is the only vertex in $e$, which is not colored by $u$ in $\chi_{i-1}$.
\end{enumerate}
	Then we try to recolor $v_i$ according to the value of the random variable $\eta_i$:
	\begin{itemize}
	\item if $\eta_i=0$, then we do not recolor $v_i$,
	\item if $\eta_i\ne 0$, then we recolor $v_i$ in the color $\eta_i$.
\end{itemize}
	In all the other situations, we do not change the color of $v_i$. Let the resulting coloring be $\chi_i$.
\end{enumerate}

Let $\tilde{\chi}$ be the coloring obtained after the consideration of all the vertices.	

\vspace{2mm}

Now we are going to give a more formal construction of the random coloring
$\tilde{\chi}$ using the techniques of random variables. This is very useful for the further proof. We analyze the event ${\cal F}$ that $\tilde{\chi}$ is not a proper coloring for $H$. We divide ${\cal F}$ into some ``local'' events and estimate their probabilities. Finally, we use Local Lemma to show that all these events do not occur simultaneously with positive probability. This implies that $\tilde{\chi}$ is a proper coloring of $H$ with positive probability, and, hence, $H$ is $r$-colorable.

\subsection{Formal Construction of the random coloring from $\S$3.2}
Without loss of generality, we may assume, that $V=\{ 1,2,3,\ldots,w\}$. Let us also fix an arbitrary ordering of the edges of $H$. Consider, on some probability space, the following set of mutually independent random elements:
\begin{itemize}
\item[1.] $\xi_1,\ldots,\xi_w$ --- equally distributed random variables, taking values $1,2,\ldots,r$ with equal probability $1/r$.

\item[2.] $\eta_1,\ldots,\eta_w$ --- equally distributed random variables taking values $1,2,\ldots,r$ with equal probability $p$ and the value $0$ with probability $1-rp$. We take the parameter $p$ equal to $p=q/(r-1)$. By the condition (\ref{condition2}) such choice of the parameter is correct, i. e., for every $r\ge 2$, one has the inequalities $rp\le r/(2(r-1))\le 1$.
\end{itemize}
\vspace{2mm}

Let $e\in E$ be an edge of $H$. For every $u=1,\ldots,r$, let ${\cal M}(e,u)$ and ${\cal AM}(e,u)$ denote the following events:
\begin{equation}\label{t1_00}
	{\cal M}(e,u)=\bigcap_{s\in\, e}\{\xi_s=u\},\quad {\cal AM}(e,u)=\left\{0<\sum_{s\in\, e}I\{\xi_s\ne u\}\le t+\omega-2\right\}.
\end{equation}
We shall introduce successively random variables $\zeta_i$, $i=1,\ldots,w$. Let ${\cal D}_1$ and ${\cal A}_1$ denote the following events:
$$
	{\cal D}_1=\bigcup_{e\in\,E:\,1\in\,e}\bigcup_{u=1}^r{\cal M}(e,u),
$$
$$
	{\cal A}_1=\bigcup_{f\in\,E:\;1\in\, f}\bigcup_{u=1}^r\left(\left\{\xi_1\ne u,\;\eta_1=u,\; \sum_{s\in\,f:\;s>1}I\{\xi_s=u\}=k-1\right\}\cap{\cal AM}(f,u)\right),
$$
and let
$$
	\zeta_1=\xi_1\,I\{\overline{{\cal D}_1}\cup\{\eta_1=0\}\cup{\cal A}_1\}+\eta_1\,I\{{\cal D}_1\cap\{\eta_1\ne 0\}\cap \overline{{\cal A}_1}\}.
$$
For every $i>1$, let ${\cal D}_i$ and ${\cal A}_i$ denote the events
$$
	{\cal D}_i=\bigcup_{e\in\,E:\,i\in\,e}\bigcup_{u=1}^r\left\{{\cal M}(e,u)\cap\bigcap_{s\in\,e:\,s<i}\{\zeta_s=u\}\right\},
$$
$$
	{\cal A}_i=\bigcup_{f\in\,E:\;i\in f}\bigcup_{u=1}^r\left(\left\{\xi_i\ne u,\;\eta_i=u,\; \sum_{s\in\,f:\;s<i}I\{\zeta_s=u\}+\sum_{s\in\,f:\;s>i}I\{\xi_s=u\}=k-1\right\}\cap{\cal AM}(f,u)\right).
$$
We define $\zeta_i$ in the following way:
$$
	\zeta_i=\xi_i\,I\{\overline{{\cal D}_i}\cup\{\eta_i=0\}\cup{\cal A}_i\}+\eta_i\,I\{{\cal D}_i\cap\{\eta_i\ne 0\}\cap \overline{{\cal A}_i}\}.
$$
It is easy to see that the random variables $\zeta_i$ take values only from $\left\{1,2,\ldots,r\right\}$. So, we may interpret the random vector $\vec{\zeta}=\left(\zeta_1,\ldots,\zeta_w\right)$ as a random $r$-coloring of the vertex set $V$ (we assign the color $\zeta_i$ to the vertex $i$). Let $\cal F$ denote the event that $\vec{\zeta}$ is not a  proper coloring of the hypergraph $H$, i. e.,
\begin{equation}\label{t1_01}
	{\cal F}=\bigcup_{e\in\, E}\bigcup_{u=1}^r\bigcap_{s\in\, e}\{\zeta_s=u\}.
\end{equation}
Our task is to prove that ${\sf P}({\cal F})<1$ under the conditions of Theorem \ref{theorem 4}.
\par
We shall divide the event $\bigcap_{s\in\, e}\{\zeta_s=u\}$ into three parts, depending on the behavior of the random variables $\{\xi_s:\;s\in e\}$. Let ${\cal C}_0(e,u)$, ${\cal C}_1(e,u)$, ${\cal C}_2(e,u)$ be the following events:
$$
	{\cal C}_0(e,u)=\bigcup_{a=1,\,a\ne u}^r\bigcap_{s\in\, e}\{\zeta_s=u,\;\xi_s=a\},\quad
	{\cal C}_1(e,u)=\bigcap_{s\in\, e}\{\zeta_s=u,\;\xi_s=u\},
$$
\begin{equation}\label{t1_02}
	{\cal C}_2(e,u)=\bigcap_{s\in\, e}\{\zeta_s=u\}\cap\bigcap_{a=1}^r\overline{{\cal M}(e,a)}.
\end{equation}
We shall consider these events separately. But before we establish a simple inequality, which we will use later. It follows from (\ref{condition2}) that
\begin{equation}\label{simple_n}
	\alpha\ln k=qk \ge 2.
\end{equation}
Note that the last inequality implies that the parameter $t$ in (\ref{t}) is correctly defined (there is no negative number under the square root).

\subsection{First part of ${\cal F}$: the event ${\cal C}_0(e,u)$}

If the event ${\cal C}_0(e,u)$ occurs, then for every $s\in e$, one has $\zeta_s=\eta_s$, since $\zeta_s\ne\xi_s$. We get the relation
\begin{equation}\label{t1_03}
	\bigcup_{u=1}^r{\cal C}_0(e,u)\subset \bigcup_{u=1}^r\bigcup_{a=1,\,a\ne u}^r\bigcap_{s\in\, e}\{\eta_s=u,\;\xi_s=a\}=
	{\cal Q}_0(e).
\end{equation}
The probability of the event ${\cal Q}_0(e)$ can be easily calculated:
\begin{equation}\label{t1_q0}
	{\sf P}({\cal Q}_0(e))=\sum_{u=1}^r\sum_{a=1,\,a\ne u}^r \prod_{s\in\, e}{\sf P}\left(\eta_s=u,\;\xi_s=a\right)=r(r-1)\left(\frac pr\right)^k.
\end{equation}

\subsection{Second part of ${\cal F}$: the event ${\cal C}_1(e,u)$}

Suppose that the event ${\cal C}_1(e,u)$ occurs. This event, obviously, implies all the events ${\cal D}_s,\; s\in e$. Then the equality $\xi_s=\zeta_s=u$ for a vertex $s\in e$ can happen in two ways: either $\eta_s\in\{0,u\}$, or $\eta_s\notin\{0,u\}$ and the event ${\cal A}_s$ occurs. Consider the following partition of the event ${\cal C}_1(e,u)$:
\begin{equation}\label{t1_04}
	{\cal C}_1(e,u)={\cal S}_0(e,u)\cup{\cal S}_1(e,u),
\end{equation}
where
$$
	{\cal S}_0(e,u)={\cal C}_1(e,u)\cap\left\{\sum_{s\in\,e}I\left\{\eta_s\notin\{0,u\}\right\}\le t+\omega-1\right\},
$$
$$
	{\cal S}_1(e,u)={\cal C}_1(e,u)\cap\left\{\sum_{s\in\,e}I\left\{\eta_s\notin\{0,u\}\right\}> t+\omega-1\right\}.
$$
Consider the event ${\cal S}_0(e,u)$. By the definition (\ref{t1_02}) of the event ${\cal C}_1(e,u)$ the following relation holds:
$$
	{\cal S}_0(e,u)\subset \bigcap_{s\in\, e}\{\xi_s=u\}\cap\left\{\sum_{s\in\,e}I\left\{\eta_s\notin\{0,u\}\right\}\le t+\omega-1\right\}.
$$
Let ${\cal Q}_1(e)$ denote the union of the last events:
\begin{equation}\label{t1_05}
	\bigcup_{u=1}^r{\cal S}_0(e,u)\subset \bigcup_{u=1}^r\left\{\bigcap_{s\in\, e}\{\xi_s=u\}\cap
	\left\{\sum_{s\in\,e}I\left\{\eta_s\notin\{0,u\}\right\}\le t+\omega-1\right\}\right\}={\cal Q}_1(e).
\end{equation}
The probability of ${\cal Q}_1(e)$ has the following estimate:
$$
	{\sf P}\left({\cal Q}_1(e)\right)=r^{1-k}\sum_{j=0}^{t+\omega-1} {k\choose j} q^{j} (1-q)^{k-j}\le r^{1-k} (1-q)^{k-t-\omega} \sum_{j=0}^{t+\omega-1} (kq)^{j}\le
$$
\begin{equation}\label{t1_q1}	
\le
	r^{1-k} (1-q)^{k-t-\omega} (kq)^{t+\omega}.
\end{equation}
The last inequality follows from the bound (\ref{simple_n}): $kq=\alpha\ln k\ge 2$.
\par
Consider now the event ${\cal S}_1(e,u)$. Let us fix $v\in e$ satisfying $\eta_v\notin\{0,u\}$. As it was noted above, the event  ${\cal A}_v$ should happen for every such vertex. This event implies that for some edge $f$, $v\in f$, and some color $a\ne u$, the following event has to occur
$$
	{\cal W}(v,f,u,a)=\left\{\xi_v=u,\;\eta_v=a,\;\sum_{s\in\,f:\;s<v}I\{\zeta_s=a\}+ \sum_{s\in\,f:\;s>v}I\{\xi_s=a\}=k-1 \right\}\cap{\cal AM}(f,a).
$$
It is easy to show that $f\ne e,$ moreover, $f\cap e = \{v\}$. Indeed, for all $s\in e$, it holds that  $\xi_s=\zeta_s=u$, but for all $s\in f\backslash \{v\}$, either $\zeta_s=a$, or $\xi_s=a$.
\par
Suppose $\{v_1,\ldots,v_h\}=\{v\in e:\; \eta_v\notin\{0,u\}\}$. For any $i=1,\ldots,h$, the event ${\cal S}_1(e,u)$ implies the event ${\cal W}(v_i,f_i,u,a_i)$ for some edge $f_i$ satisfying $\{v_i\}=f_i\cap e$ and some color $a_i\ne u$. Moreover, ${\cal S}_1(e,u)$ also implies that $h=h(e,u)\ge t+\omega$. Since there are at most $\omega$ 3-cycles containing $e$, there is a subset $\{f'_1,\ldots,f'_t\}\subset \{f_1,\ldots,f_h\}$ such that $f'_i$ and $f'_j$ are disjoint for all $i\ne j$.

For further convenience, we introduce a notation of \emph{the configuration of the first type}. For given edge $e$, the set of edges $\{f_1,\ldots,f_t\}$ is said to be the configuration of the first type (denotation: $\{f_1,\ldots,f_t\}\in {\rm CONF1}(e)$) if, for any any $i=1,\ldots,t$, $|f_i\cap e|=1$ and, moreover, all the edges $f_i$ are pairwise disjoint.

Thus, by the above arguments we the following relation
\begin{equation}\label{t1_06}
	{\cal S}_1(e,u)\subset \bigcap_{s\in\, e}\{\xi_s=u\}\cap\bigcup_{\substack{a_1,\ldots,\,a_t=1\\a_i\ne u}}^r\bigcup_{\{f_1,\ldots,\,f_t\}\in\, {\rm CONF1}(e)}\bigcap_{i=1}^t
	{\cal W}(e\cap f_i,f_i,u,a_i),
\end{equation}
where the set of edges $\{f_1,\ldots,f_t\}$ is assumed to be ordered according to the originally selected ordering of $E$, i. e. the number of the edge $f_i$ is less than the number of the edge $f_j$, if $i<j$. Let us use the notations: $\widehat{f}_i=f_i\backslash e$ and $v_i=e\cap f_i$, $i=1,\ldots,t$. It follows from the definition of the configuration of the first type that the sets $\widehat{f}_i$, $i=1,\ldots,t$ do not have common vertices, i.e. $\widehat{f}_i\cap \widehat{f}_j=\emptyset$, if $i\ne j$. Furthermore, $|\widehat{f}_i|=k-1$.

If the event ${\cal W}(e\cap f_i,f_i,u,a_i)$ happens, then by ${\cal AM}(f_i,a_i)$ the edge $f_i$ contains at most $t+\omega-2$ vertices $s$, satisfying $\xi_s\ne a_i$. Moreover, for all such vertices, $\zeta_s=a_i$, and so, $\zeta_s=\eta_s=a_i$. The set $\widehat{f}_i$ contains at most $t+\omega-3$ such vertices, since the vertex $v_i$ doesn't belong to $\widehat{f}_i$ and $\xi_{v_i}=u\ne a_i$. Thus, we obtain the relation
$$
	\bigcap_{s\in\, e}\{\xi_s=u\}\cap\bigcap_{i=1}^t{\cal W}(e\cap f_i,f_i,u,a_i)\subset \bigcap_{s\in\, e}\{\xi_s=u\}\cap\bigcap_{i=1}^t\{\eta_{v_i}=a_i\}\cap
$$
\begin{equation}\label{t1_07}
	\cap\bigcap_{i=1}^t\left\{
	\bigcap_{s\in\,\widehat{f}_i}\left(\{\xi_s\ne a_i,\;\eta_s=a_i\}\cup\{\xi_s=a_i\}\right) \right\}\cap\bigcap_{i=1}^t\left\{\sum_{s\in\,\widehat{f}_i}I\{\xi_s\ne a_i\}\le t+\omega-3\right\}.
\end{equation}
Let ${\cal Q}_2(e,F)$ denote the union of the last events, where $F=\{f_1,\ldots,f_t\}\in {\rm CONF1}(e)$:
$$
{\cal Q}_2(e,F)=\bigcup_{u=1}^r\bigcup_{\substack{a_1,\ldots,\,a_t=1\\a_i\ne u}}^r\left\{\bigcap_{s\in\, e}\{\xi_s=u\}\cap\bigcap_{i=1}^t\{\eta_{v_i}=a_i\}\cap\right.
$$
\begin{equation}\label{t1_08}	
	\left.\cap\bigcap_{i=1}^t\left\{
	\bigcap_{s\in\,\widehat{f}_i}\left(\{\xi_s\ne a_i,\;\eta_s=a_i\}\cup\{\xi_s=a_i\}\right) \right\}\cap
	\bigcap_{i=1}^t\left\{\sum_{s\in\,\widehat{f}_i}I\{\xi_s\ne a_i\}\le t+\omega-3\right\}\right\}.
\end{equation}
The relations (\ref{t1_06}) and (\ref{t1_07}) imply
\begin{equation}\label{t1_09}
\bigcup_{u=1}^r {\cal S}_1(e,u)\subset \bigcup_{F\in\; {\rm CONF1}(e)}{\cal Q}_2(e,F).
\end{equation}
Let us estimate the probability of ${\cal Q}_2(e,F)$:
$$
	{\sf P}\left({\cal Q}_2(e,F)\right)=\sum_{u=1}^r\sum_{\substack{a_1,\ldots,\,a_t=1\\a_i\ne u}}^r r^{-k} p^t \prod_{i=1}^t\sum_{j=0}^{t+\omega-3} {\vert \widehat{f}_i\vert\choose j} \left(\frac {r-1}r\right)^jp^j \left(\frac 1r\right)^{\vert \widehat{f}_i\vert-j}=
$$
$$
	=r(r-1)^t r^{-k} p^t r^{-\sum\limits_{i=1}^t\vert \widehat{f}_i\vert}\prod_{i=1}^t\sum_{j=0}^{t+\omega-3} {\vert \widehat{f}_i\vert\choose j} q^j=
	r(r-1)^t r^{-k} p^t r^{-t(k-1)}\prod_{i=1}^t\sum_{j=0}^{t+\omega-3} {k-1\choose j} q^j
	\le
$$	
\begin{equation}\label{t1_q2}		
	\le r^{-(t+1)(k-1)} q^t \prod_{i=1}^t\sum_{j=0}^{t+\omega-3} k^j q^j
	\le r^{-(t+1)(k-1)} q^t (kq)^{t(t+\omega-2)}.
\end{equation}

\subsection{Third part of ${\cal F}$: the event ${\cal C}_2(e,u)$}

We shall show that if the event ${\cal C}_2(e,u)$ happens then the sum $\sum_{s\in e}I\{\xi_s\ne u\}$ cannot be very small. We shall establish the equality
\begin{equation}\label{t1_11}
{\cal C}_2(e,u)={\cal C}_2(e,u)\cap\left\{\sum_{s\in\, e}I\{\xi_s\ne u\}\ge t+\omega-1\right\}.
\end{equation}
Indeed, let us consider the intersection of three events (see the definition of the event ${\cal C}_2(e,u)$ in (\ref{t1_02})):
$$
	{\cal C}_2(e,u)\cap\left\{\sum_{s\in\, e} I\{\xi_s\ne u\}\le t+\omega-2\right\}=
$$
$$
=\bigcap_{s\in\, e}\{\zeta_s=u\}\cap\bigcap_{a=1}^r\overline{{\cal M}(e,a)}\cap\left\{\sum_{s\in\, e} I\{\xi_s\ne u\}\le t+\omega-2\right\}.
$$
The second and the third events imply the happening of the event ${\cal AM}(e,u)$ (see (\ref{t1_00})). The first one implies that for every $s\in e$ satisfying $\xi_s\ne u$, we have $\zeta_s=\eta_s=u$. Moreover, since the event ${\cal AM}(e,u)$ holds, the set of such vertices in not empty. Consider a vertex $v\in e$ satisfying $\xi_v\ne u$ and $\xi_s=u$ for every $s\in e,\; s>v$. It is clear that the event ${\cal A}_v$ holds. So, $\zeta_v=\xi_v\ne u$, and we get a contradiction with the first event in the intersection. Thus, these three events are inconsistent, and we prove the equality (\ref{t1_11}).
\par
Let us estimate the probability of ${\cal C}_2(e,u)$. Consider the random set $T=\{s\in e:\; \xi_s\ne u\}$. The event ${\cal C}_2(e,u)$ implies, first, that all $v\in T$ satisfy $\zeta_v=\eta_v=u$, and second, that $\vert T\vert\ge t+\omega-1$ (see (\ref{t1_11})). Let us use the denotation: $E(e)=\{f\in E\backslash \{e\}:\; f\cap e\neq\emptyset\}$.

If $\zeta_v\ne\xi_v$ for some vertex $v$, then there should happen at least two events: the event $D_v$ and the event
\begin{equation}
	{\cal B}(e,f_v,v,u,a_v)=\left\{{\cal M}(f_v,a_v)\cap\bigcap_{s\in\,f_v:\,s<v}\{\zeta_s=a_v\} \cap\{\zeta_v=\eta_v=u\} \right\},
\label{t1_12}	
\end{equation}
where $f_v$ is some edge, satisfying $v\in e\cap f_v$, $v$ is the first vertex from $e\cap f_v$ and $a_v\ne u$ is some color. It is clear that  edges $f_v$ are different for different $v$.
\par
Let $Y$ be an arbitrary subset of the edge $e$ satisfying $y=\vert Y\vert \ge t+\omega-1$.  Then we have the inclusion
\begin{equation}\label{t1_13}
	{\cal C}_2(e,u)\cap\{T=Y\}\subset \bigcap_{s\in\, e\backslash Y}\{\xi_s=u\}
	\cap\bigcup_{\substack{a_1,\ldots,\,a_{y}=1\\a_i\ne u}}^r\bigcup_{\substack{f_1,\ldots,f_{y}\in\,E(e)\\v_i\,\in\,f_i\cap e}}\bigcap_{i=1}^{y}{\cal B}(e,f_i,v_i,u,a_i),
\end{equation}
where $v_i$ is the first vertex from $f_i\cap e$. Since $y\ge t+\omega-1$ and the edge $e$ is contained in at most $\omega$ 3-cycles, there is a set of edges $\{f'_1,\ldots,f'_{t-1}\}\subset\{f_1,\ldots,f_y\}$ such that the sets $\widehat{f'_i}=f'_i\setminus e$, $i=1\ldots,t-1$, are pairwise disjoint and, moreover, the first vertices of $f_i\cap e$ are different for different $i=1,\ldots,t-1$.

For further convenience, we introduce a notation of \emph{the configuration of the second type}. For given edge $e$, an unordered set of edges $F=\{f_1,\ldots,f_{t-1}\}$ is said to be the configuration of the second type (denotation: $F\in {\rm CONF2}(e)$) if, for any any $i=1,\ldots,t$, $f_i\in E(e)$ and, moreover, all the sets $f_i\setminus e$ are pairwise disjoint.
Let us also use the following denotation:
$$
  S(F)=\bigcup_{i=1}^{t-1}\left(f_i\cap e\right),
$$
where $F=\{f_1,\ldots,f_{t-1}\}\in {\rm CONF2}(e)$.

Due to (\ref{t1_13}) and the above argument, we have the following inclusion:
$$
	{\cal C}_2(e,u)\subset \bigcup_{\substack{F\in\;{\rm CONF}2(e)\\F=\{f_1,\ldots,f_{t-1}\}}}
\left\{\bigcap_{s\in\, e\backslash S(F)}\left\{\{\xi_s=u\}\sqcup\{\xi_s\ne u,\;\eta_s=u\}\right\}\cap\right.
$$
\begin{equation}\label{t1_14}
	\left.\cap\bigcap_{s\in\, S(F)}\{\eta_s=u\} \cap\bigcup_{\substack{a_1,\ldots,\,a_{t-1}=1\\a_i\ne u}}^r\bigcap_{i=1}^{t-1}{\cal B}(e,f_i,v_i,u,a_i)\right\},
\end{equation}
where the set of edges $\{f_1,\ldots,f_t\}$ is assumed to be ordered according to the originally selected ordering of $E$, and, moreover, for any $i=1,\ldots,t-1$, $v_i$ denotes the first vertex in $f_i\cap e$.

\par
The event ${\cal B}(e,f_i,v_i,u,a_i)$ is obviously (see (\ref{t1_12})) contained in the event $\bigcap_{s\in f_i}\{\xi_s=a_i\}\cap\{\eta_{v_i}=u\}$. Hence, by (\ref{t1_14}) we get the relation
$$
	{\cal C}_2(e,u)\subset \bigcup_{\substack{F\in\;{\rm CONF}2(e)\\F=\{f_1,\ldots,f_{t-1}\}}}
	\left\{\bigcap_{s\in\, e\backslash S(F)}\left\{\{\xi_s=u\}\sqcup\{\xi_s\ne u,\;\eta_s=u\}\right\}\cap\right.
$$
$$
	\left.\cap\bigcap_{s\in\, S(F)}\{\xi_s\ne u,\;\eta_s=u\}\cap\bigcup_{\substack{a_1,\ldots,\,a_{t-1}=1\\a_i\ne u}}^r\bigcap_{i=1}^{t-1}\bigcap_{s\in f_i}\{\xi_s=a_i\}\right\}.
$$
Let us take the union of both parts over $u$:
$$
	\bigcup_{u=1}^r{\cal C}_2(e,u)\subset \bigcup_{\substack{F\in\;{\rm CONF}2(e)\\F=\{f_1,\ldots,f_{t-1}\}}}\bigcup_{u=1}^r
	\left\{\bigcap_{s\in\, e\backslash S(F)}\left\{\{\xi_s=u\}\sqcup\{\xi_s\ne u,\;\eta_s=u\}\right\}\cap\right.
$$
\begin{equation}\label{t1_16}
	\left.\cap\bigcap_{s\in\, S(F)}\{\eta_s=u\}\cap\bigcup_{\substack{a_1,\ldots,\,a_{t-1}=1\\a_i\ne u}}^r\bigcap_{i=1}^{t-1}\bigcap_{s\in f_i}\{\xi_s=a_i\}\right\}.
\end{equation}
\par
Let us introduce the following event:
$$
	{\cal Q}_3(e,F)=\bigcup_{u=1}^r\left\{\bigcap_{s\in\, e\backslash S(F)}\left\{\{\xi_s=u\}\sqcup\{\xi_s\ne u,\;\eta_s=u\}\right\}\cap\right.
$$
\begin{equation}\label{t1_18}
	\left.\cap\bigcap_{s\in\, S(F)}\{\eta_s=u\}\cap\bigcup_{\substack{a_1,\ldots,\,a_{t-1}=1\\a_i\ne u}}^r\bigcap_{i=1}^{t-1}\bigcap_{s\in f_i}\{\xi_s=a_i\}\right\},
\end{equation}
where $e$ is an edge of $H$, $F=\{f_1,\ldots,f_{t-1}\}\in {\rm CONF}2(e)$, and the edges are written according to the original ordering. By (\ref{t1_16}) and (\ref{t1_18}) we get the relation
\begin{equation}\label{t1_19}
\bigcup_{u=1}^r{\cal C}_2(e,u)\subset \bigcup_{F\in\;{\rm CONF}2(e)} {\cal Q}_3(e,F).
\end{equation}
Now we are going to estimate the probability of ${\cal Q}_3(e,F)$.

\par
Let us consider more closely the set of edges $F=\{f_1,\ldots,f_{t-1}\}$. The hypergraph $H(F)=(V,F)$ can be divided into some number of connected components. Suppose $H_1,\ldots,H_l$ are these components. Since $F$ is a configuration of the second type, the edges $f_i$ and $f_j$ can have common vertices only inside the edge $e$. Moreover, we know that $H$ is 2-simple. For every component $H_j$, $j=1,\ldots,l$, let us use the denotations:
$$
  h_j=|\{f\in E(H_j):\; |f\cap e|=2\}|\mbox{ and }l_j=|\{f\in H_j:\; |f\cap e|=1\}|.
$$
Due to the 2-simplicity of $H$ we have
\begin{equation}\label{t1_new1}
  \sum_{j=1}^l(h_j+l_j)=t-1.
\end{equation}
For any component $H_j$, let $G_j$ be the following graph: $G_j=(V_j,E_j)$, where
$$
  V_j=e\cap V(H_j),\quad E_j=\{f\cap e:\; f\in E(H_j)\mbox{ and }|f\cap e|=2\}|.
$$
The following claim clarifies the structure of the configurations of the second type.

\vspace{2mm}
\begin{claim}\label{claim_tree} For any $j=1,\ldots,l$, $G_j$ is either a tree or an isolated vertex.
\end{claim}
\begin{proof1} If $h_j=0$ then $l_j>0$. But in this case $H_j$ consists of only one edge $f$. Indeed, if $g\in E(H_j)$, $g\ne f$, then $|g\cap e|=1$ and, moreover, $|g\cap f|>0$. This can only happens when $f\cap e=g\cap e$, i.e. $g$ has the same first vertex in its intersection with $e$ as $f$. This fact is in conflict with the definition of the configuration of the second type. So, $G_j$ is just an isolated vertex.

Now let $h_j>0$. Since $H_j$ is connected and $F\in {\rm CONF}2(e)$, $G_j$ is also connected.
Suppose there is a cycle $(w_1,\ldots,w_m)$, $m\ge 3$, in $G_j$, i.e. $\{w_i,w_{i+1}\}\in E(G_j)$, $i=1,\ldots,m-1$ and, moreover, $\{w_1,w_{m}\}\in E(G_j)$. Without loss of generality, assume that $w_1<w_j$ for any $j>1$, i.e. $w_1$ is the vertex with the least number in the cycle. Since $\{w_1,w_2\}\in E(G_j)$, there is an edge $g_1\in E(H_j)$ such that $\{w_1,w_2\}=g_1\cap e$, so $w_1$ is the first vertex in $g_1\cap e$. By analogy, there is another edge $g_2\in E(H_j)$ such that $\{w_1,w_m\}=g_2\cap e$, so $w_1$ is also the first vertex in $g_2\cap e$. We obtain a contradiction with the fact that $F\in {\rm CONF}2(e)$. Hence, $G_j$ is a tree.
\end{proof1}

Claim \ref{claim_tree} implies that $|V(G_j)|=|E(G_j)|+1=h_j+1$, thus,
\begin{equation}\label{t1_new2}
  |S(F)|=\sum_{j=1}^l |V(G_j)|=\sum_{j=1}^l(h_j+1).
\end{equation}
Moreover, from the definitions of the values $h_j$ and $l_j$ we get that, for any $j=1,\ldots,l$,
\begin{equation}\label{t1_new3}
  \left|\bigcup_{f\in H_j}f\right|=(k-2)h_j+(k-1)l_j+|V(G_j)|=(k-1)(h_j+l_j)+1.
\end{equation}
Finally, Claim \ref{claim_tree} implies that, for any $j=1,\ldots,l$,
\begin{equation}\label{t1_new4}
  l_j\le 1.
\end{equation}
Indeed, since there is a bijection between the edges of $H_j$ and the first vertices in their intersection with $e$, we have
$$
  h_j+l_j=|E(H_j)|\le |V(G_j)|=h_j+1,
$$
and the inequality (\ref{t1_new4}) follows.

Using the notations introduced above one can easily find the probability of the event ${\cal Q}_3(e,F)$ (see (\ref{t1_18})):
\begin{equation}\label{t1_21}
	{\sf P}\left({\cal Q}_3(e,F)\right)= r\left(\frac 1r+\frac qr\right)^{k-|S(F)|}\,p^{|S(F)|}
  (r-1)^l\prod_{j=1}^l r^{-\left|\bigcup_{f\in H_j}f\right|}.
\end{equation}
Let us explain the last two factors in the right-hand side of (\ref{t1_21}). Since all the edges of $F$ are monochromatic in the main coloring $\xi$, the values of $\xi_{s}$ should coincide for all $s\in V(H_j)$. Thus, we only have to choose a color (not equal to $u$) for every component (the factor $(r-1)^l$). The last factor is equal to the probability (we have already chosen the colors) that every edge in the component $H_j$ is monochromatic in the main coloring $\xi$.

Using obtained relations (\ref{t1_new1}), (\ref{t1_new2}), (\ref{t1_new3}), (\ref{t1_new4}), we get the following estimate of the probability of the event ${\cal Q}_3(e,F)$:
$$
	{\sf P}\left({\cal Q}_3(e,F)\right)= r\left(\frac 1r+\frac qr\right)^{k-|S(F)|}\,p^{|S(F)|}
  (r-1)^l\prod_{j=1}^l r^{-\left|\bigcup_{f\in H_j}f\right|}=
$$
$$
  =r^{1-k}(1+q)^{k-|S(F)|}(rp)^{|S(F)|}(r-1)^l \prod_{j=1}^l r^{-(k-1)(h_j+l_j)-1}=
$$
$$
  = r^{1-k}(1+q)^{k-|S(F)|}(rp)^{|S(F)|}(r-1)^l r^{-(k-1)(t-1)-l}\le
$$
\begin{equation}\label{t1_q3}
  \le r^{(1-k)t}(1+q)^k (rp)^{|S(F)|}\le r^{(1-k)t}(1+q)^k (2q)^{t-1}.
\end{equation}
We need to comment only the last inequality. From the condition (\ref{condition2}) we have $2q\le 1$ and, moreover, we know that $rp=(r/(r-1))q\le 2q$. Finally, from (\ref{t1_new1}), (\ref{t1_new3}) and (\ref{t1_new4}) we immediately see that $|S(F)|\ge t-1$.

The bound (\ref{t1_q3}) completes the estimation of different parts of the event ${\cal F}$. Now we shall prove that the probability of ${\cal F}$ is less than 1 under the conditions of Theorem \ref{theorem 4}.

\subsection{Application of Local Lemma for estimating the probability of ${\cal F}$}

Remember that by the definitions (\ref{t1_01}) and (\ref{t1_02}) of the events ${\cal F}$ and ${\cal C}_i(e,u)$, $i=1,2,3$, $e\in E$, $u=1,\ldots,r$, we have the equality
$$
	{\cal F}=\bigcup_{e\in E}\bigcup_{u=1}^r\left({\cal C}_1(e,u)\cup{\cal C}_2(e,u)\cup{\cal C}_3(e,u)\right).
$$
It follows from the obtained relations (\ref{t1_03}), (\ref{t1_04}), (\ref{t1_05}), (\ref{t1_09}) and (\ref{t1_19}), that
\begin{equation}\label{t1_22}
	{\cal F}\subset\bigcup_{e\in\, E}\left\{{\cal Q}_0(e)\cup{\cal Q}_1(e)\right\}\cup
	\bigcup_{e\in\, E}\bigcup_{F\in\,{\rm CONF}1(e)}{\cal Q}_2(e,F)\cup
	\bigcup_{e\in\, E}\bigcup_{F\in\,{\rm CONF}2(e)}{\cal Q}_3(e,F).
\end{equation}

Further, we shall use a classical claim, which is called Local Lemma. This statement was first proved in the paper of
P. Erd\H{o}s and L. Lov\'asz \cite{ErdLov}. We shall formulate it in a special case.

\vspace{2mm}
\begin{theorem}\label{local_lemma} Let events ${\cal B}_1,\ldots,$ ${\cal B}_N$ be given on some probability space. Let $S_1,\ldots,S_N$ be subsets of ${\cal R}_N=\{1,\ldots,N\}$ such that for any $i=1,\ldots,N$, the event ${\cal B}_i$ is independent of the algebra generated by the events $\{{\cal B}_j,j\in {\cal R}_N\backslash S_i\}$. If, for any $i=1,\ldots,N$, the following inequality holds
\begin{equation}\label{loclemma}
  \sum\limits_{j\in\, S_i}{\sf P}({\cal B}_j)\le 1/4,
\end{equation}
then ${\sf P}\left(\bigcap\limits_{j=1}^N \overline{{\cal B}_j}\right)\ge \prod\limits_{j=1}^N\left(1-2{\sf P}({\cal B}_i)\right)>0$.
\end{theorem}
\vspace{2mm}

The proof of the Local Lemma can be found in the book \cite{AlonSpencer}. Consider the system of events $\Psi$ consisting of all the events ${\cal Q}_i(e)$, $i=0,1$, $e\in E$, the events ${\cal Q}_2(e,F)$, $e\in E$, $F\in\,{\rm CONF}1(e)$, and also the events ${\cal Q}_3(e,F)$, $e\in E$, $F\in\,{\rm CONF}2(e)$. By (\ref{t1_22}) the inequality holds
\begin{equation}\label{t1_23}
	{\sf P}\left({\cal F}\right)\le {\sf P}\left(\bigcup_{{\cal B}\in\,\Psi}{\cal B}\right)=1-{\sf P}\left(\bigcap_{{\cal B}\in\,\Psi}\overline{{\cal B}}\right).
\end{equation}
We shall show that the probability of $\bigcap_{{\cal B}\in\,\Psi}\overline{{\cal B}}$ is greater than zero. Due to Local Lemma (see Theorem \ref{local_lemma}), it is sufficient to find, for every ${\cal B}\in\Psi$, a system of events $\Psi_{\cal B}\subset\Psi$ such that ${\cal B}$ and the algebra generated by $\left\{{\cal Q}\in\Psi\backslash \Psi_{\cal B}\right\}$ are independent, and, moreover, such that the following inequality holds:
\begin{equation}\label{t1_24}
    \sum\limits_{{\cal Q}\in\,\Psi_{\cal B}}{\sf P}({\cal Q})\le 1/4.
\end{equation}

\bigskip
\par
The event ${\cal B}\in\Psi$ can have one of the following three types:
\begin{enumerate}
  \item ${\cal B}={\cal Q}_i(e)$ for some $e\in E$ and $i\in\{0,1\}$;
  \item ${\cal B}={\cal Q}_2(e,F)$ for some $e\in E$ and $F\in\,{\rm CONF}1(e)$;
  \item ${\cal B}={\cal Q}_3(e,F)$ for some $e\in E$ and $F\in\,{\rm CONF}2(e)$.
\end{enumerate}
For any ${\cal B}\in\Psi$, we define \emph{the domain} $D({\cal B})$ of the event ${\cal B}$ as follows:
$$
  D({\cal B})=\begin{cases} e,& \text{if }{\cal B}={\cal Q}_i(e),\;i=0,1;\\
  e\cup \left(\bigcup\limits_{f\in\, F}f\right),& \text{if }{\cal B}={\cal Q}_i(e,F),\;i=2,3.
  \end{cases}
$$
By the definitions (\ref{t1_03}), (\ref{t1_05}), (\ref{t1_08}), (\ref{t1_18}) the event ${\cal B}$ belongs to the algebra generated by the random variables $\{\xi_j,\eta_j:\;j\in D({\cal B})\}$. Then this event is independent of the algebra generated by the random variables $\{\xi_j,\eta_j:\;j\in V\setminus D({\cal B})\}$. We take the system $\Psi_{\cal B}$ consisting of all the events ${\cal Q}\in\Psi$ such that the domains of ${\cal Q}$ and ${\cal B}$ have nonempty intersection. In other words,
$$
  \Psi_{\cal B}=\left\{{\cal Q}\in\Psi:\;\;D({\cal Q})\cap D({\cal B})\ne\emptyset\right\}.
$$
Thus, the event ${\cal B}$ is independent of the algebra generated by $\left\{{\cal J}\in\Psi\backslash \Psi_{\cal B}\right\}$, if we choose $\Psi_{\cal B}$ in this way. We have to check the inequality (\ref{t1_24}). By the choice of the set $\Psi_{\cal B}$ we get the relation
$$
	\sum\limits_{{\cal J}\in\,\Psi_{\cal B}}{\sf P}\left({\cal J}\right)\le\sum_{e\in E:\;e\cap\,D({\cal B})\ne\emptyset }\left({\sf P}\left({\cal Q}_0(e)\right)+{\sf P}\left({\cal Q}_1(e)\right)\right)+\sum_{\substack{e\,\in E,\;F\in\,{\rm CONF}1(e):\\D({\cal B})\cap D({\cal Q}_2(e,F))\ne\emptyset}} {\sf P}\left({\cal Q}_2(e,F)\right)+
$$
\begin{equation}\label{t1_25}
+\sum_{\substack{e\,\in E,\;F\in\,{\rm CONF}2(e):\\D({\cal B})\cap D({\cal Q}_3(e,F))\ne\emptyset}} {\sf P}\left({\cal Q}_3(e,F)\right).
\end{equation}
Let us denote by $a({\cal B})$, $b({\cal B})$ and $c({\cal B})$ the number of summands in the first sum, the second sum and the third sum in the right-hand side of (\ref{t1_25}) respectively. Using these denotations from the relation (\ref{t1_25}) and the estimates (\ref{t1_q0}), (\ref{t1_q1}), (\ref{t1_q2}), (\ref{t1_q3}) we get the inequality
$$
	\sum\limits_{{\cal J}\in\,\Psi_{\cal B}}{\sf P}({\cal B})\le a({\cal B})\left(r(r-1)\left(\frac pr\right)^{k}+r^{1-k}(1-q)^{k-t-\omega}(kq)^{t+\omega}\right)+
$$
\begin{equation}\label{t1_26}
	+b({\cal B})r^{-(t+1)(k-1)}q^t(kq)^{t(t+\omega-2)}+c({\cal B})r^{-t(k-1)}(1+q)^k(2q)^{t-1}.
\end{equation}

Now we shall consider three cases depending on ${\cal B}$.

\begin{enumerate}
  \item ${\cal B}={\cal Q}_i(e)$ for some $e\in E$ and $i\in\{0,1\}$. By the condition (\ref{d_def}) of Theorem \ref{theorem 4} there exist at most $d$ other edges intersecting an arbitrary $e\in E$. So,
      $$
        a({\cal B})\le d+1,\;b({\cal B})\le (d+1){d\choose t}+(d+1)d{d-1\choose t-1},
      $$
      \begin{equation}\label{t1_27}
        c({\cal B})\le (d+1){d\choose t-1}+(d+1)d{d-1\choose t-2}.
      \end{equation}
      The first inequality in (\ref{t1_27}) is obvious. To show the last two it is sufficient to notice that $e$ can intersect either with $e'$ from the event ${\cal Q}_2(e',F)$ or with some $f\in F$.

  \item ${\cal B}={\cal Q}_2(e,F)$ for some $e\in E$ and $F\in\,{\rm CONF}1(e)$. This event depends on $(t+1)$ edges. So, by using the estimates from (\ref{t1_27}) we get
       $$
        a({\cal B})\le (t+1)(d+1),\;b({\cal B})\le (t+1)\left((d+1){d\choose t}+(d+1)d{d-1\choose t-1}\right),
      $$
      \begin{equation}\label{t1_27_1}
        c({\cal B})\le (t+1)\left((d+1){d\choose t-1}+(d+1)d{d-1\choose t-2}\right).
      \end{equation}
  \item ${\cal B}={\cal Q}_3(e,F)$ for some $e\in E$ and $F\in\,{\rm CONF}2(e)$. This event depends on $t$ edges. So, as in the previous case
       $$
        a({\cal B})\le t(d+1),\;b({\cal B})\le t\left((d+1){d\choose t}+(d+1)d{d-1\choose t-1}\right),
      $$
      \begin{equation}\label{t1_27_2}
        c({\cal B})\le t\left((d+1){d\choose t-1}+(d+1)d{d-1\choose t-2}\right).
      \end{equation}
\end{enumerate}

It is easy to see from (\ref{t1_27}), (\ref{t1_27_1}) and (\ref{t1_27_2}) that the maximal upper bounds for $a({\cal B})$, $b({\cal B})$ and $c({\cal B})$ are in the second case. So, to prove (\ref{t1_24}) it is sufficient to establish (due to (\ref{t1_26})) the following inequality:
$$
 W=(t+1)(d+1)\left(r(r-1)\left(\frac pr\right)^{k}+r^{1-k}(1-q)^{k-t-\omega}(kq)^{t+\omega}\right)+
$$
$$
  +(t+1)\left((d+1){d\choose t}+(d+1)d{d-1\choose t-1}\right)r^{-(t+1)(k-1)}q^t(kq)^{t(t+\omega-2)}+
$$
\begin{equation}\label{t1_n2}
	+(t+1)\left((d+1){d\choose t-1}+(d+1)d{d-1\choose t-2}\right)r^{-t(k-1)}(1+q)^k(2q)^{t-1}\le 1/4
\end{equation}

We shall need some additional estimates contained in the next section.

\subsection{Auxiliary analytics}
The value $W$ (see (\ref{t1_n2})) consists of four summands:
$$
	(t+1)(d+1)r(r-1)\left(p/r\right)^k,\quad
	(t+1)(d+1)r^{1-k} (1-q)^{k-t-\omega} (kq)^{t+\omega},
$$
$$
	(t+1)\left((d+1){d\choose t}+(d+1)d{d-1\choose t-1}\right)q^t r^{-(k-1)(t+1)}(kq)^{t(t+\omega-2)},
$$
$$
	(t+1)\left((d+1){d\choose t-1}+(d+1)d{d-1\choose t-2}\right)r^{-(k-1)t} (1+q)^k (2q)^{t-1}.
$$
Consider and estimate them separately.
\begin{enumerate}

\item The first summand is $(t+1)(d+1)r(r-1)\left(p/r\right)^k$. Using the restriction (\ref{d_def}), the conditions (\ref{condition1}) and (\ref{condition2}) we obtain the upper bound for the first summand:
$$
(t+1)(d+1)r(r-1)\left(\frac pr\right)^n\le (t+1)kr^{k-1}r^{1-n}(r-1) \left(\frac q{r-1}\right)^k=
$$
\begin{equation}\label{ineq_1}
=(t+1)k (r-1)^{1-k} q^k\le k^2 q^k\le k^2 2^{-k}.
\end{equation}
	
\item The second summand is $(t+1)(d+1)r^{1-k} (1-q)^{k-t-\omega} (kq)^{t+\omega}$. Since the choice of parameter $q$ in (\ref{t}), we get the relations
$$
	(t+1)(d+1)r^{1-k} (1-q)^{k-t-\omega} (kq)^{t+\omega}\le (t+1)k r^{k-1}r^{1-k} (1-q)^{k-t-\omega} (kq)^{t+\omega}=
$$
$$
	= (t+1)n(1-q)^{k-t-\omega} \left(\alpha\ln k\right)^{t+\omega}\le (t+1) k e^{q(t+\omega)-qk} \left(\alpha\ln k\right)^{t+\omega}=
$$
\begin{equation}\label{ineq_2}
	=(t+1)k^{1-\alpha} e^{\alpha(\ln k)(t+\omega)/k}\left(\alpha\ln k\right)^{t+\omega}.
\end{equation}	
	
\item Let us consider the third summand in the expression (\ref{t1_n2}) for the value $W$:
\begin{equation}\label{t1_29}
	(t+1)\left((d+1){d\choose t}+(d+1)d{d-1\choose t-1}\right)q^t r^{-(k-1)(t+1)}(kq)^{t(t+\omega-2)}.
\end{equation}
We shall need some preliminary estimates.
\par
First, the following inequalities hold:
$$
	(t+1)\left((d+1){d\choose t}+(d+1)d{d-1\choose t-1}\right)=(t+1){d\choose t}(d+1)(t+1)\le
$$
\begin{equation}\label{t1_30}
\le (t+1)^2(d+1)\frac {d^t}{t!}\le (t+1)^2\frac {(d+1)^t}{t!}.
\end{equation}
Second, the choice of parameters $t$ and $q$ (see (\ref{t})) implies the relations
$$
	q^t (kq)^{t(t+\omega-2)}=k^{-t}(kq)^{t(t+\omega-1)}\le k^{-t}(kq)^{t^2+t\omega}= k^{-t} \exp\left\{t^2\ln \left(\alpha\ln k\right)\right\}(kq)^{t\omega}\le
$$
\begin{equation}\label{t1_32}
\le k^{-t}\exp\left\{\ln k\right\}(\alpha\ln k)^{t\omega} = k^{1-t}(\alpha\ln k)^{t\omega}.
\end{equation}
Finally, from (\ref{t1_30}), (\ref{t1_32}) and the original restriction (\ref{d_def}), we obtain the upper bound for the expression (\ref{t1_29}):
$$
	(t+1)\left((d+1){d\choose t}+(d+1)d{d-1\choose t-1}\right)q^t r^{-(k-1)(t+1)}(kq)^{t(t+\omega-2)}\le
$$
$$
	\le \frac {(t+1)^2}{t!} (d+1)^{t+1} r^{-(k-1)(t+1)} k^{1-t}(\alpha\ln k)^{t\omega}\le
	\frac {(t+1)^2}{t!} k^{(t+1)(1-b/t)} k^{1-t}(\alpha\ln k)^{t\omega}\le
$$
$$
  \le \frac {(t+1)^2}{t!}\, k^{t+1-(b(t+1)/t)} k^{1-t}(\alpha\ln k)^{t\omega}=
$$
\begin{equation}\label{ineq_3}	
	=\frac {(t+1)^2}{t!}\, k^{2-b-(b/t)}(\alpha\ln k)^{t\omega}\le \frac {(t+1)^2}{t!}\, k^{2-b}(\alpha\ln k)^{t\omega}.
\end{equation}

\item It remains to estimate the fourth summand in the expression (\ref{t1_n2}) for the value $W$:
\begin{equation}\label{t1_33}
(t+1)\left((d+1){d\choose t-1}+(d+1)d{d-1\choose t-2}\right)r^{-(k-1)t}(1+q)^k(2q)^{t-1}.
\end{equation}
By an analogy with (\ref{t1_30}), we get:
$$
(t+1)\left((d+1){d\choose t-1}+(d+1)d{d-1\choose t-2}\right)=(t+1){d\choose t-1}(d+1)t\le
$$
\begin{equation}\label{t1_34}
\le (t+1)t\left(\frac {de}{t-1}\right)^{t-1}(d+1)\le (d+1)^t \left(\frac {e}{t-1}\right)^{t-1} (t+1)t.
\end{equation}
Further, by (\ref{t}) it holds that
\begin{equation}\label{t1_35}
	(2q)^{t-1}(1+q)^{k}\le k^{1-t} (2\alpha\ln k)^{t-1}e^{qk}=k^{1+\alpha-t} (2\alpha\ln k)^{t-1}.
\end{equation}
Finally, from (\ref{t1_34}), (\ref{t1_35}) and (\ref{d_def}) we obtain an upper bound for the expression (\ref{t1_33}):
$$
	(t+1)\left((d+1){d\choose t-1}+(d+1)d{d-1\choose t-2}\right)r^{-(k-1)t}(1+q)^k(2q)^{t-1}\le
$$
$$
	\le (t+1)t \left(\frac {e}{t-1}\right)^{t-1} (d+1)^t r^{-(k-1)t} k^{1+\alpha-t}(2\alpha\ln k)^{t-1}\le
$$
$$
	\le (t+1)t \left(\frac {2e\alpha\ln k}{t-1}\right)^{t-1} r^{(k-1)t}\,k^{t(1-b/t)} r^{-(k-1)t}k^{1+\alpha-t}=
$$
\begin{equation}\label{ineq_4}
	=(t+1)t \left(\frac {2e\alpha\ln k}{t-1}\right)^{t-1} k^{1+\alpha-b}.
\end{equation}

The inequality (\ref{ineq_4}) completes the estimation of the parts of the value $W$.

\end{enumerate}

\subsection{The completion of the proof of Theorem \ref{theorem 4}}

Let us gather the obtained bounds for the summands in the expression (\ref{t1_n2}) for the value $W$. The relations (\ref{ineq_1}), (\ref{ineq_2}), (\ref{ineq_3}) and (\ref{ineq_4}) imply the inequalities
$$
	W\le \frac {k^2}{2^k}+(t+1)k^{1-\alpha}e^{\alpha(\ln k)(t+\omega)/k}\left(\alpha\ln k\right)^{t+\omega}+\frac {(t+1)^2}{t!}\, k^{2-b}(\alpha\ln k)^{t\omega}+
$$
$$
+(t+1)t\left(\frac {2e\alpha\ln k}{t-1}\right)^{t-1} k^{1+\alpha-b}<\frac 14,
$$
the last of which holds, since the condition (\ref{condition3}) of theorem \ref{theorem 4}. Thus, the required relation  (\ref{t1_n2}) is established. It implies the inequality (\ref{t1_24}) necessary for the application of Local Lemma. It follows from Local Lemma that the probability of simultaneous happening of all the events $\overline{{\cal B}}$, where ${\cal B}\in\Psi$, is greater than zero. Then by (\ref{t1_23}) we have shown that
$$
	{\sf P}\left({\cal F}\right)<1.
$$
Let us complete the proof. Indeed, we have proved, that the probability of the event that the random coloring $\vec{\zeta}$ is not a proper coloring of $H$ is less than one. So, $\vec{\zeta}$ is a proper coloring with positive probability, and $\chi(H)\le r$. Theorem \ref{theorem 4} is proved.

\subsection{The completion of the proof of Theorem \ref{theorem 3}}

We shall use Theorem \ref{theorem 4}. Let us choose the parameters $b$ and $\alpha$:
$$
	b=4,\quad \alpha=2.
$$
By this choice of $b$, $\alpha$ and the condition $\omega\le\sqrt{\ln k/(\ln\ln k)}$ there exists an integer $k_1$ such that for all $k\ge k_1$, the inequalities (\ref{condition1}) and (\ref{condition2}) hold. Let us consider the left part of (\ref{condition3}). We have $t=O\left(\sqrt{\ln k/\ln\ln k}\right)$ (see (\ref{t})), so
$$
	(t+1)k^{1-\alpha} e^{\alpha(\ln k) (t+\omega)/k}  \left(\alpha\ln k\right)^{t+\omega}=e^{O\left(\ln\ln k\right)}k^{-1}e^{o(1)}e^{O\left(\sqrt{\ln k\ln \ln k}\right)}=o(1),\; k\to\infty,
$$
$$
	\frac {(t+1)^2}{t!}\, k^{2-b}(\alpha\ln k)^{t\omega}=O\left(k^{-2}\right)e^{\ln k(1+o(1))}=o(1),\; k\to\infty,
$$
$$
	(t+1)t\left(\frac {2e\alpha\ln k}{t-1}\right)^{t-1} k^{1+\alpha-b}=e^{O\left(\sqrt{\ln k\ln\ln k}\right)} k^{-1}=o(1),\; k\to\infty.
$$
These relations imply the existence of an integer $k_2$ such that the inequality (\ref{condition3}) holds, for all $k\ge k_2$.

Let $H=(V,E)$ be an $k$-uniform hypergraph, $H\in\mathcal{H}(k,r,\omega)$ with $\omega\le\sqrt{\ln k/(\ln\ln k)}$. In the case $k\ge k_0=\max(k_1,k_2)$ the hypergraph $H$ satisfies all the conditions of Theorem \ref{theorem 4}, except (\ref{d_def}). But $H$ is not $r$-colorable, and so there exists an edge $e\in E$ with edge degree at least $\left\lfloor r^{k-1}k^{1-b/t}\right\rfloor$. So, the edge $e$ contains a vertex with degree at least
$$
	\left\lfloor r^{k-1}k^{1-b/t}\right\rfloor/k +1\ge \left(r^{k-1}k^{1-b/t}-1\right)/k +1=
	r^{k-1}k^{-b/t}+1-1/k.
$$
Thus, we have established the inequality $\Delta(H)>r^{k-1}k^{-b/t}$ and, consequently,
$$
	\Delta(\mathcal{H}(n,r,\omega))\ge r^{k-1}k^{-b/t}=r^{k-1}k^{-4\left\lfloor\sqrt{\frac{\ln k}{\ln(2\ln k)}}\right\rfloor^{-1}}.
$$
Theorem \ref{theorem 3} is proved.

\section{Choosability in random hypergraphs}

In this section we will discuss $r$-choosability of the random hypergraph $H(n,k,p)$. Let us recall the required definitions.

We say that a hypergraph $H$ is $r$-\textit{choosable} if for every family of sets $L=\{L(v):\; v\in V\}$ ($L$ is called \textit{list assignment}), such that $\vert L(v)\vert=r$ for all $v\in V$, there is a proper coloring from the lists (for every $v\in V$ we should use a color from $L(v)$). The \textit{choice number} of a hypergraph $H$, denoted by $ch(H)$, is the least $r$ such that $H$ is $r$-choosable. It is clear that $\chi(H)\le ch(H)$. The choice numbers of graphs were independently introduced by V.G. Vizing (see \cite{Vizing}) and by P. Erdos, A. Rubin and H. Taylor (see \cite{ERT}). In this paper we consider the threshold probability for $r$-choosability of $H(n,k,p)$.

\subsection{Threshold for $r$-choosability in $H(n,k,p)$}

The choice number of the random hypergraph $H(n,k,p)$ has been studied by M. Krivelevich and V. Vu (see \cite{KrivVu}). They proved that $ch(H(n,k,p))$ is asymptotically very closed to $\chi(H(n,k,p))$. Their first result holds for almost all $p$, but has at little gap between $ch(H(n,k,p))$ and $\chi(H(n,k,p))$.

\begin{theorem}\label{thm_krivvu_1}{\rm (M. Krivelevich, V. Vu, \cite{KrivVu})} Suppose $k\ge 2$ is fixed. There exists a constant $C=C(k)$ such that, for any $p$ satisfying $C\,n^{1-k}\le p\le 0.9$, the following convergence holds
$$
  {\sf P}\left(ch(H(n,k,p))\le (1+\psi(n))\,k^{1/(k-1)}\,\chi(H(n,k,p))\right)\to 1\mbox { as }n\to\infty,
$$
where $\psi(n)\to 0$ as $n^{k-1}p\to\infty$.
\end{theorem}

The second theorem from \cite{KrivVu} states that, for sufficiently large $p$, this gap can be removed.

\begin{theorem}\label{thm_krivvu_2}{\rm (M. Krivelevich, V. Vu, \cite{KrivVu})} Suppose $k\ge 2$ is fixed and $0<\varepsilon<(k-1)^2/(2k)$. Then, for any $p$ satisfying $n^{-(k-1)^2/(2k)+\varepsilon}\le p\le 0.9$, the following convergence holds
$$
  {\sf P}\left(ch(H(n,k,p))=(1+o(1))\,\chi(H(n,k,p))\right)\to 1\mbox { as }n\to\infty.
$$
\end{theorem}

Theorem \ref{thm_krivvu_1} together with Theorem \ref{thm_krivsud} of Krivelevich and Sudakov implies the following corollary, which is an analogue of Corollary \ref{cor_krivsud} for $r$-choosability.

\begin{corollary}\label{cor_krivvu} Let $k\ge 3$ and $\varepsilon\in(0,1)$ be fixed. There is a constant $r_0=r_0(k,\varepsilon)$ such that, for any $r=r(n)$ satisfying the conditions
$$
  r\ge r_0,\quad r^{k-1}\ln r=o(n^{k-1}),
$$
the following convergence holds:
$$
  {\sf P}\left(ch(H(n,k,p))\le r\right)\to 1,\mbox{ where }p=(1-\varepsilon)\,\frac{r^{k-1}\ln r}k\,\frac{n}{{n\choose k}}.
$$
\end{corollary}
We see that the lower bound for the threshold for $r$-choosability provided by Corollary \ref{cor_krivvu} does not coincide with the upper bound in Lemma \ref{lemma_2} (if hypergraph is not $r$-colorable then it is also not $r$-choosable). Their ratio has an order of $k$. Recall that applying Theorem \ref{thm_krivsud} provides the following restrictions on the parameters $r$ and $k$ in Corollary \ref{cor_krivvu} (see (\ref{cor_ks_condition})):
\begin{equation}\label{cor_kv_condition_1}
  r=\Omega\left(k^{29}(\ln k)^{28}\right),\quad n\ge k^{9k+O(k\ln\ln k/(\ln k))}.
\end{equation}
What can be said about the lower bound for the threshold for $r$-choosability when $r=O\left(k^{29}(\ln k)^{28}\right)$?

\bigskip
\begin{note} It should be noted that, for very large $r$ (e.g., $r>\sqrt{n}$) and fixed $k$, Theorem \ref{thm_krivvu_2} together with Theorem \ref{thm_krivsud} gives an asymptotic value for the required threshold for $r$-choosability:
$$
  p^*\sim r^{k-1}\ln r\,\frac{n}{{n\choose k}}.
$$
\end{note}

\bigskip
The proof of Theorem \ref{thm_akkt_rcolors} by Achlioptas, Krivelevich, Kim and Tetali is based on the determi\-nistic coloring algorithm, which cannot be generalized to the case of an arbitrary $r$-uniform list assignment, so the lower bound (\ref{akkt_rcolors}) does not hold for the threshold probability for $r$-choosability. Moreover, the proof of the result (\ref{ach_moore}) by Achlioptas and Moore is also cannot be adopted for list colorings. Thus, in the case when $r$ is small in comparison with $k$ we have only the result of Lemma \ref{lemma_1} which is just a generalization of the result (\ref{alon_spencer_01}) by Alon and Spencer.

\begin{lemma}\label{lemma_4} There exists an integer $k_0$ and a positive number $c$ such that for any fixed $k\ge k_0$ and $r\ge 2$, the following statement holds:
\begin{equation}\label{sta_lemma_4}
\mbox{ if }p\le c\,\frac{r^{k-1}}{k^2}\,\frac {n}{{n\choose k}},\mbox{ then }{\sf P}\left(H(n,k,p)\mbox { is $r$-choosable}\right)\to 1.
\end{equation}
\end{lemma}

Our new results concerning $r$-choosability in random hypergraphs are formulated in the following two theorems.

\begin{theorem}\label{thm_list_1} Suppose $k=k(n)\ge 3$ and $r=r(n)\ge 3$ satisfy the relation
$$
  \frac{3}{128}\,\frac{r^{k-1}}{\sqrt{k}}-\ln n\to -\infty\mbox{ as }n\to\infty.
$$
If
\begin{equation}\label{thresh_list_bound1}
  p\le \frac 3{32}\,\frac{r^{k-1}}{k^{3/2}}\,\frac n{{n\choose k}},
\end{equation}
then ${\sf P}\left(H(n,k,p)\mbox { is $r$-choosable}\right)\to 1$.
\end{theorem}

\vspace{2mm}
\begin{theorem}\label{thm_list_2}
Suppose $\delta\in(0,1)$ is a constant. Let $k=k(n)$ and $r=r(n)\ge 2$ satisfy the following conditions: $k\ge k_0$ , where $k_0$ is some absolute constant, and, moreover,
$$
(k-1)\ln r < \frac {1-\delta}2\ln n,\quad r^{k-1}k^{-\varphi(k)}\ge 6\ln n,
$$
where $\varphi(k)=4\left\lfloor\sqrt{\frac{\ln k}{\ln(2\ln k)}}\right\rfloor^{-1}$. Then for  function $p=p(n)$, satisfying
\begin{equation}\label{thresh_list_bound2}
p\le \frac 12\,\frac {r^{k-1}}{k^{1+\varphi(k)}}\,\frac n{{n\choose k}},
\end{equation}
we have ${\sf P}\left(H(n,k,p)\mbox{ is $r$-choosable}\right)\to 1$ as $n\to\infty$.
\end{theorem}

It is easy to see that Theorem \ref{thm_list_1} and Theorem \ref{thm_list_2} stated that the results of Corollary \ref{corollary1} (assertion 1)) and Theorem \ref{theorem 2} also hold in the case of list colorings. Both provided bounds are better (for all sufficiently large $k$) than the result (\ref{sta_lemma_4}) of Lemma \ref{lemma_4}, and both are worse than the result of Corollary \ref{cor_krivvu}. Hence, Theorem \ref{thm_list_2} gives the best lower bound (\ref{thresh_bound2}) for the threshold probability for $r$-choosability of $H(n,k,p)$ in the wide area of the parameters (recall the restriction (\ref{cor_kv_condition_1})):
$$
  r\le k^{29}(\ln k)^{28}\;\mbox{ and }\; 6\,k^{\varphi(k)}\ln n\le r^{k-1}\le n^{(1-\delta)/2},
$$
where $k\ge k_0$ is sufficiently large. The inequality (\ref{thresh_bound1}), in comparison with
(\ref{thresh_bound2}), does not have an upper restriction $r^{k-1}\le n^{(1-\delta)/4}$, so it provides the best lower bound in the area
$$
  3\le r\le k^{29}(\ln k)^{28}\;\mbox{ and }\; r^{k-1}\ge n^{(1-\delta)/2}.
$$

The proofs of Theorems \ref{thm_list_1} and Theorem \ref{thm_list_2} are very similar to the proofs of the first assertion of Corollary \ref{corollary1} and Theorem \ref{theorem 2}, so we do not give the complete argument and describe only the main ides and differences.

\subsection{Ideas of the proofs of Theorems \ref{thm_list_1} and \ref{thm_list_2}}

For given $k,r\ge 2$, let $\Delta_{list}(k,r)$ denote the minimum possible $\Delta(H)$, where $H$ is a $k$-uniform non-$r$-choosable hypergraph. In \cite{Shab1} D.A. Shabanov shows that the lower bound (\ref{d(k,r)_shab}) for $\Delta(k,r)$ (see $\S$2.1) holds for $\Delta_{list}(k,r)$ also: for any $k,r\ge 3$,
$$
  \Delta_{list}(k,r)>\frac 18\,k^{-1/2}r^{k-1}.
$$
Using this inequality one can easily prove Theorem \ref{thm_list_1} by the same argument as in Lemma \ref{lemma_3}.

\bigskip
\begin{note} The lower bound (\ref{d(k,r)_kkr}) for $\Delta(k,r)$ obtained by Kostochka, Kumbhat and R\"odl does not hold for $\Delta_{list}(k,r)$, so we cannot apply it to $r$-choosability of random hypergraphs.
\end{note}

\bigskip
To prove Theorem \ref{thm_list_2} it is sufficient to show that under the conditions of Theorem \ref{theorem 4} the hypergraph $H$ is not only $r$-colorable, but is $r$-choosable. The proof of $r$-choosability remains almost the same. The difference appears in the distributions of the random variables.

Suppose $H=(V,E)$ is a $k$-uniform hypergraph satisfying the conditions of Theorem \ref{theorem 4} and let $L=\{L(v):\; v\in V\}$ be an $r$-uniform list assignment with the set of colors ${\mathbb N}$. Without loss of generality, $V=\{1,\ldots,w\}$. In comparison with $\S$3.3 we introduce random variables with another distribution. Let $\xi_1,\ldots,\xi_w$ and $\eta_1,\ldots,\eta_w,$ be mutually independent random variables with the following distribution:
\begin{itemize}
	\item $\xi_i$, $i=1,\ldots,w$, has the uniform distribution on the set $L(i)$ ($i=1,\ldots,w$),
	
	\item $\eta_i$, $i=1,\ldots,w$, takes all values from $L(i)$ with the same probability $p$ and the value $0$ with probability $1-rp$.
\end{itemize}

For given edge $e\in E$, let $M(e)$ be equal to $M(e)=\bigcap_{s\in\, e}L(s)$. For every $u\in M(e)$, we introduce the events ${\cal M}(e,u)$, ${\cal AM}(e,u)$ whose definitions are the same as in $\S$3.3 (see (\ref{t1_00})). Then we construct the random coloring $\vec{\zeta}=\left(\zeta_1,\ldots,\zeta_w\right)$ by the same way as in $\S$3.3. The only difference is that in the definitions of the events ${\cal A}_i$ and ${\cal D}_i$ the parameter $u$ does not take values from 1 to $r$, it should take values from an appropriate set $M(e)$ or $M(f)$. The rest of the proof remains the same without any unobvious change.

\renewcommand{\refname}{References}


\begin{thebibliography}{99}

\bibitem{BolTom} B. Bollob\'as, A. Thomason, ``Thresholds functions'', \emph{Combinatorica}, \textbf{7} (1987), 35--38.

\bibitem{Friedgut} E. Friedgut, ``Necesary and sufficient conditions for sharp thresholds of graph properties'', \emph{Journal of the American Mathematical Society}, \textbf{12} (1999), 1017--1054.

\bibitem{AlonSpencerUnp} N. Alon, J. Spencer, ``A note on coloring random $k$-sets'', Unpublished manuscript.

\bibitem{AKKT} D. Achlioptas, J.H. Kim, M. Krivelevich, P. Tetali, ``Two-colorings random hypergraphs'', \emph{Random Structures and Algorithms}, \textbf{20}:2(2002), 249--259.

\bibitem{AchMoore} D. Achlioptas, C. Moore, ``On the 2-colorability of random hypergraphs'', \emph{Random}, Springer-Verlag, Berlin, 2002, 78--90.

\bibitem{KarLuc} M. Karo\'nski, T. $\L$uczak, ``Random hypergraphs'', \emph{Combinatorics, Paul Erd\H{o}s is eighty}, Bolyai Society Mathematical Studies, \textbf{2}, D. Mikl\'os, V.T. S\'os, T. Sz\H{o}nyi eds., 1996, 283--293.

\bibitem{KrivSud} M. Krivelevich, B. Sudakov, ``The chromatic numbers of random hypergraphs'', \emph{Random Structures and Algorithms}, \textbf{12} (1998), 381--403.

\bibitem{ErdLov} P. Erd\H{o}s, L. Lov\'asz, ``Problems and results on 3-chromatic hypergraphs and some related questions'', \emph{Infinite and Finite Sets}, Colloquia Mathematica Societatis Janos Bolyai, \textbf{10}, North Holland, Amsterdam, 1973, 609--627.

\bibitem{KostRodl} A. V. Kostochka, V. R\"odl, ``Constructions of sparse uniform hypergraphs with high chromatic number'', \emph{Random Structures and Algorithms}, \textbf{36}:1(2010), 46--56.

\bibitem{RadhSrin} J. Radhakrishnan, A. Srinivasan, ``Improved bounds and algorithms for hypergraph two-coloring'', \emph{Random Structures and Algorithms}, \textbf{16}:1(2000), 4--32.

\bibitem{Shab1} D. A. Shabanov, ``On $r$-chromatic hypergraphs'', \emph{Discrete Mathematics}, to appear.

\bibitem{KKR} A.V. Kostochka, M. Kumbhat, V. R\"odl, ``Coloring uniform hypergraphs with small edge degrees'', \emph{Fete of Combinatorics and Computer Science}, Bolyai Society Mathematical Studies, \textbf{20}, Springer, 2010, 213--238.

\bibitem{JLR} S. Jansen, T. $\L$uczak, A. Rucinski, \emph{Random graphs}, Wiley-Interscience, New York, 2000.

\bibitem{Szabo} Z. Szab\'o, ``An application of Lovasz Local Lemma - a new lower bound for the van der Waerden number'', \emph{Random Structures and Algorithms}, \textbf{1}:3(1990), 343--360.

\bibitem{KostKumb} A. V. Kostochka, M. Kumbhat, ``Coloring uniform hypergraphs with few edges'', \emph{Random Structures and Algorithms}, \textbf{35}:3(2009), 348--368.

\bibitem{Shab2} D. A. Shabanov, ``Random coloring method in the combinatorial problem of Erd\H{o}s and Lov\'asz'', \emph{Random Structures and Algorithms}, to appear.

\bibitem{Beck} J. Beck, ``On 3-chromatic hypergraphs'', \emph{Discrete Mathematics}, \textbf{24}:2(1978), 127--137.

\bibitem{Spencer} J. H. Spencer, ``Coloring n-sets red and blue'', \emph{J. Combinatorial Theory, Series A}, \textbf{30}(1981), 112--113.

\bibitem{AlonSpencer} N. Alon, J. H. Spencer, \emph{Probabilistic method}, Wiley--Interscience, New York, 2002.

\bibitem{Vizing} V. G. Vizing, ``Coloring the vertices of a graph in prescribed colors'', \emph{Metody diskretnogo analiza v teorii kodov i skhem} (in Russian), Sobolev Institute of Mathematics, Novosibirsk, \textbf{29}(1976), 3--10.

\bibitem{ERT} P. Erd\H{o}s, A. L. Rubin, H. Taylor, ``Choosability in graphs'', \emph{Proc. West Coast Conference on Combinatorics, Graph Theory and Computing}, \textbf{26}, 1980, 125--157.

\bibitem{KrivVu} M. Krivelevich, V. Vu, ``Choosability in random hypergraphs'', \emph{Journal of Combinatorial Theory, Series B}, \textbf{83}:2(2001), 241--257.


\end{thebibliography}
\end{document}